\documentclass[journal]{IEEEtran}
\usepackage{array}
\usepackage{marvosym}
\usepackage[cmex10]{amsmath}
\usepackage{amssymb,amsthm,epsfig}
\usepackage{graphpap}
\usepackage[latin1]{inputenc}
\usepackage{rotating}
\usepackage{algorithm}
\usepackage{setspace}
\usepackage{blkarray}
\usepackage{environ}
\usepackage{multirow}
\usepackage{setspace}
\usepackage{textcomp}
\usepackage[flushleft]{threeparttable}
\usepackage{color}
\usepackage{soul}

\usepackage{epsfig}
\usepackage{amsmath}
\usepackage{algcompatible}
\usepackage{kbordermatrix}
\usepackage{latexsym}
\usepackage{amsfonts}
\usepackage[latin1]{inputenc}
\usepackage{rotating}
\usepackage{multirow}
\usepackage{setspace}
\usepackage{textcomp}
\usepackage{bm}
\hypersetup{colorlinks=true, allcolors=black, linkcolor=black, citecolor=black }
\usepackage{mathtools}
\DeclarePairedDelimiter{\ceil}{\lceil}{\rceil}
\DeclarePairedDelimiter{\floor}{\lfloor}{\rfloor}
\usepackage{mdwlist}
\usepackage{graphicx}
\usepackage{cite}
\usepackage{enumitem}

\usepackage{footnote}
\makesavenoteenv{tabular}
\usepackage{tabulary}
\usepackage[table]{xcolor}
\newcolumntype{P}[1]{>{\centering\arraybackslash}p{#1}}
\newcolumntype{M}[1]{>{\centering\arraybackslash}m{#1}}
\usepackage{booktabs, makecell}
\usepackage[export]{adjustbox}
\pdfoutput=1

\DeclareMathOperator*{\argmax}{arg\,max}
\DeclareMathOperator*{\argmin}{arg\,min}
\begin{document}
\title{Optimizing Service Restoration in Distribution Systems with Uncertain Repair Time and Demand}
%
\author{Anmar Arif,~\IEEEmembership{Student Member,~IEEE,}
		Shanshan Ma,~\IEEEmembership{Student Member,~IEEE,}
       Zhaoyu Wang,~\IEEEmembership{Member,~IEEE,}
       \\Jianhui Wang,~\IEEEmembership{Senior Member,~IEEE,}
       Sarah M. Ryan,~\IEEEmembership{Senior Member,~IEEE,}       
       Chen Chen,~\IEEEmembership{Member,~IEEE,}

\thanks{This work was partially supported by the U.S. Department of Energy Office of Electricity Delivery and Energy Reliability, the National Science Foundation under grant ECCS1609080.}
\thanks{A. Arif is with the Department of Electrical and Computer Engineering, Iowa State University, Ames, IA, 50011 USA, and also with the Department of Electrical Engineering, King Saud University, Riyadh, 11451 Saudi Arabia (Email:aiarif\MVAt iastate.edu).}
\thanks{S. Ma, and Z. Wang are with the Department of Electrical and Computer Engineering, Iowa State University, Ames, IA, 50011 USA. (Email:aiarif\MVAt iastate.edu,sma\MVAt iastate.edu,wzy\MVAt iastate.edu).}
\thanks{S. M. Ryan is with the Department of Industrial and Manufacturing Systems Engineering, Iowa State University, Ames, IA, 50010 USA. (Email:smryan\MVAt iastate.edu).}
\thanks{J. Wang and C. Chen are with the Energy Systems Division, Argonne National Laboratory, Lemont, IL 60439 USA (Email: jianhui.wang\MVAt anl.gov,morningchen\MVAt anl.gov)}
}
\maketitle
\begin{abstract}
This paper proposes a novel method to co-optimize distribution system operation and repair crew routing for outage restoration after extreme weather events. A two-stage stochastic mixed integer linear program is developed. The first stage is to dispatch the repair crews to the damaged components. The second stage is distribution system restoration using distributed generators, and reconfiguration. We consider demand uncertainty in terms of a truncated normal forecast error distribution, and model the uncertainty of the repair time using a lognormal distribution. A new decomposition approach, combined with the Progressive Hedging algorithm, is developed for solving large-scale outage management problems in an effective and timely manner. The proposed method is validated on modified IEEE 34- and 8500-bus distribution test systems.
\end{abstract}


\section*{Nomenclature}
{
\addcontentsline{toc}{section}{Nomenclature}
\begin{description}[style=multiline,leftmargin=3cm]
\item[\textbf{Sets and Indices}]
\end{description}
\begin{description}[style=multiline,leftmargin=2.2cm] 
\item[$N$] Set of damaged components and the depot
\item[$m/n$] Indices for damaged components and the depot
\item[$c$] Index for crews
\item[$i/j$] Indices for buses
\item[$\Omega_{B}$] Set of buses
\item[$\Omega_{K(.,i)}$] Set of lines with bus $i$ as the to bus
\item[$\Omega_{K(i,.)}$] Set of lines with bus $i$ as the from bus
\item[$\Omega_{K(l)}$] Set of lines in loop $l$
\item[$\Omega_{SB}$] Set of substations
\item[$\Omega_{SW}$] Set of lines with switches
\item[$k$] Index for distribution line
\item[$t$] Index for time
\item[$\mathcal{S}$] Set of scenarios
\item[$s$] Index for scenario
\item[\textbf{Parameters}] 
\item[$e_{i,t,s}$] Active power load forecast error for load at bus $i$ and time $t$ in scenario $s$
{\color{black}\item[$M$] A sufficiently large positive number}
\item[$\mathcal{C}$] Number of crews
\item[$o_c/d_c$] Start/end point of crew $c$
\item[$P^{B_{max}}_k/Q^{B_{max}}_k$] Active/reactive power limit of line $k$ \vspace{0.05cm}
\item[$P^{G_{max}}_i/Q^{G_{max}}_i$] Active/reactive power limits of DGs
\item[$P^D_{i,t,s}/Q^D_{i,t,s}$] Diversified active/reactive demand at bus $i$ and time $t$ in scenario $s$
\item[$P^U_{i,t,s}/Q^U_{i,t,s}$] Undiversified active/reactive demand at bus $i$ and time $t$ in scenario $s$
\item[$\mathcal{T}_{m,s}$] The time needed to repair damaged component $m$ in scenario $s$
\item[$R_{k}/X_{k}$] Resistance/reactance of line $k$
\item[$T^R_{m,n}$] Travel time between $m$ and $n$
\item[$\omega_i$] Priority weight of load at bus $i$
\item[$\lambda$] The number of time steps a load needs to return to normal condition after restoration.
\end{description}
\begin{description}[style=multiline,leftmargin=3cm]
\item[\textbf{Decision Variables}]
\end{description}
\begin{description}[style=multiline,leftmargin=2cm] 
\item[$x_{m,n,c}$] Binary variable indicating whether crew $c$ moves from damaged component $m$ to $n$
\item[$\alpha_{m,c,s}$] Arrival time of crew $c$ at damaged component $m$ in scenario $s$
{\color{black}\item[$\beta^s_{i,j,t}$] Binary variable equals 1 if $i$ is the parent bus of $j$ and 0 otherwise in scenario $s$}
\item[$f_{m,t,s}$] Binary variable equal to 1 if damaged component $m$ is repaired at time $t$ in scenario $s$
{\color{black}\item[$P_{i,t,s}^L/Q_{i,t,s}^L$] Active/reactive load supplied at bus $i$ and time $t$ in scenario $s$
\item[$P_{i,t,s}^G/Q_{i,t,s}^G$] Active/reactive power generated by DG at bus $i$ in scenario $s$
\item[$P^B_{k,t,s}/Q^B_{k,t,s}$] Active/reactive power flowing on line $k$ }
\item[$u_{k,t,s}$] Binary variables indicating the status of the line $k$ at time $t$ in scenario $s$
\item[$V_{i,t,s}$] Voltage at bus $i$ and time $t$ in scenario $s$
\item[$y_{i,t,s}$] Connection status of the load at bus $i$ and time $t$ in scenario $s$
\item[$z_m$] Binary variable equal to 1 if damaged component $m$ is a critical component to repair
\end{description}

}
\section{Introduction}

\IEEEPARstart{N}ATURAL catastrophes have highlighted the vulnerability of the electric grids. In 2017, Hurricane Harvey and Hurricane Irma caused electric outages to nearly 300,000 \cite{harvey} and 15 million customers \cite{irma}, respectively. The loss of electricity after a hurricane or any natural disaster can cause significant inconvenience and is potentially life threatening. Improving outage management and accelerating service restoration are critical tasks for utilities. A crucial responsibility for the utilities is to dispatch repair crews and manage the network to restore service for customers. Relying on utility operators' experience to dispatch repair crews during outages may not lead to an optimal outage management plan. Therefore, there is a need to design an integrated framework to optimally coordinate repair and restoration.

Some research has been conducted to integrate repair and restoration in power transmission systems. In \cite{pre-hurricane}, a deterministic mixed integer linear programming (MILP) model was solved to assign repair crews to damaged components without considering the travel time. Reference \cite{DP_repair} presented a dynamic programming model for routing repair crews. Routing repair crews in transmission systems has been discussed by Van Hentenryck and Coffrin in \cite{Hent2015}. The authors presented a deterministic two-stage approach to decouple the routing and restoration models. The first stage solved a restoration ordering problem using MILP. The ordering problem formulation assumed  that only one damaged component can be repaired at each time step. The goal of the first stage was to find an optimal sequence of repairs to maximize the restored loads. The second-stage routing problem was formulated as a constraint programming model and solved using Neighborhood Search algorithms and Randomized Adaptive Decomposition.

In previous work, we developed a cluster-first route-second approach to solve the deterministic repair and restoration problem \cite{Arif2016}. However, a major challenge in solving the distribution system repair and restoration problem (DSRRP) is its stochastic nature. Predicting the repair time accurately for each damaged component is almost impossible. In this paper, we consider the uncertainty of the repair time and the customer load demand. We propose a two-stage stochastic mixed-integer program (SMIP) to solve the stochastic DSRRP (S-DSRRP). The first stage in the stochastic program is to determine the routes for each crew. The second stage models the operation of the distribution system, which includes distributed generation (DG) dispatch and network reconfiguration by controlling line switches. The routing problem is modeled as a vehicle routing problem (VRP), which has a long history in operations research \cite{laporte_2009}. The routing problem is an NP-hard combinatorial optimization problem with exponential computation time. Adding uncertainty and combining distribution system operation constraints with the routing problem further increase the complexity. To solve the large-scale S-DSRRP efficiently, a new decomposition algorithm is developed and combined with the Progressive Hedging (PH) algorithm. Our algorithm decomposes the S-DSRRP into two stochastic subproblems. The goal of the first subproblem is to find a set of damaged components that, if repaired, will maximize the served load. In the second subproblem, the repair crews are dispatched to the selected damaged components by solving S-DSRRP. The two subproblems are solved repeatedly, using parallel PH, until crews have been dispatched to repair all damaged components. The algorithm for solving the decomposed S-DSRRP is referred to as D-PH. The key contributions of this paper include: 1) improving our previously developed deterministic DSRRP formulation in \cite{Arif2016} by considering cold load pickup, and reducing the number of decision variables by refining crew routing constraints; 2) modeling the uncertainty of the repair time and the demand in DSRRP; 3) formulating a two-stage stochastic problem for repair and restoration; and 4) developing a new decomposition algorithm combined with parallel PH for solving large-scale S-DSRRP.

The rest of the paper is organized as follows. Section \ref{chap:2} states the modeling assumptions and presents the uncertainty in the model. Section \ref{chap:3} develops the mathematical formulation. In Section \ref{chap:4}, the proposed algorithm is presented. The simulation and results are presented in Section \ref{chap:5}, and Section \ref{chap:6} concludes this paper.
\vspace{-0.08cm}

\section{Modeling assumptions and uncertainty}\label{chap:2}

After a disastrous event that results in damages to the electric grid infrastructure, utilities first need to conduct damage assessment before mobilizing repair crews. Damage assessors patrol the network to locate and evaluate the damages to the grid, before the repair crews are dispatched. Damage assessment can be performed with the help of fault/outage identification algorithms, reports from customers, and aerial survey after extreme conditions. This paper is concerned with the phase after damage assessment; i.e., repairs and DG/switch operation. Hence, we assume that the locations of the damages are known from the assessment phase. Furthermore, it is assumed that the DGs in the system are controllable ones that are installed as back-up generators. In addition, each crew has the resources required to repair the damages. After determining the locations of damaged components, repair crews are dispatched to the damaged components to repair and restore the system.

In this paper, the uncertainties of repair time and load are represented by a finite set of discrete scenarios, which are obtained by sampling. The lognormal distribution is used to model the repair time, as recommended in \cite{zapata}.
%
Load uncertainty is modeled in terms of load forecast error \cite{Lu2013}. {\color{black}Define $P_{i,t}^F$ as the load forecast for the load at bus $i$ at time $t$, Fig. \ref{LoadForecast} shows an example of a 24-hour load profile. A load forecast error is generated independently for every hour. The forecast error for the load at bus $i$ and time $t$ in scenario $s$ is a realization of a truncated normal random variable $e_{i,t,s}$, so that the error is bounded using a fixed percentage (e.g., 15\%). The active demand for the load at bus $i$ and time $t$ in scenario $s$ is then obtained as follows: 
\begin{equation}
P^D_{i,t,s} = P_{i,t}^F (1+e_{i,t,s})
\end{equation} 
where a similar equation is used to obtain the corresponding realization for reactive power. By bounding the error to $\pm$15\%, equation (1) states that the actual load is within 15\% of the forecasted load. Fig. \ref{LoadScenarios} shows an example of 30 generated scenarios for one load, where $P_{i,t}^F$ is the load forecast, and $P_{i,t,s}^D$ is the generated scenario.

\begin{figure}[h!]
\vspace{-0.2cm}
\setlength{\abovecaptionskip}{0pt} 
\setlength{\belowcaptionskip}{0pt} 
  \centering
    \includegraphics[width=0.49\textwidth]{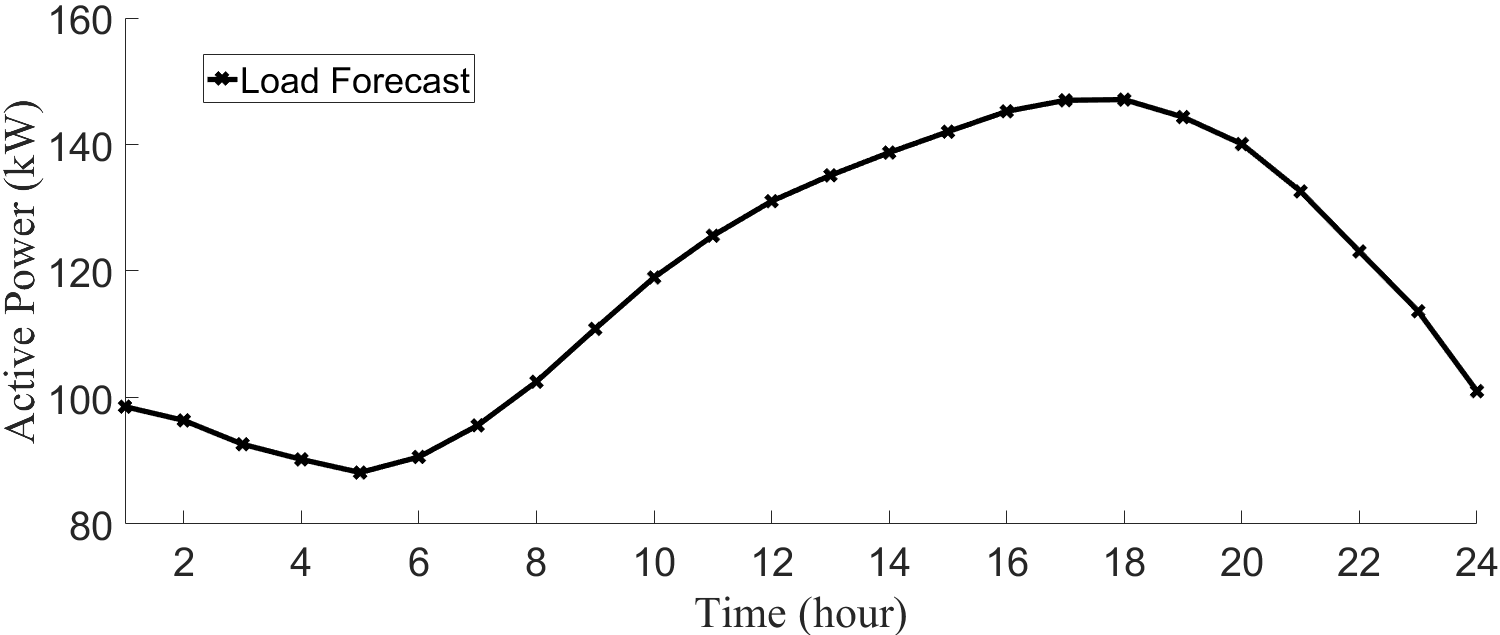}
    \caption{Forecast of active power consumption of a load.}\label{LoadForecast}
\end{figure}
\begin{figure}[h!]
\vspace{-0.2cm}
\setlength{\abovecaptionskip}{0pt} 
\setlength{\belowcaptionskip}{0pt} 
  \centering
    \includegraphics[width=0.49\textwidth]{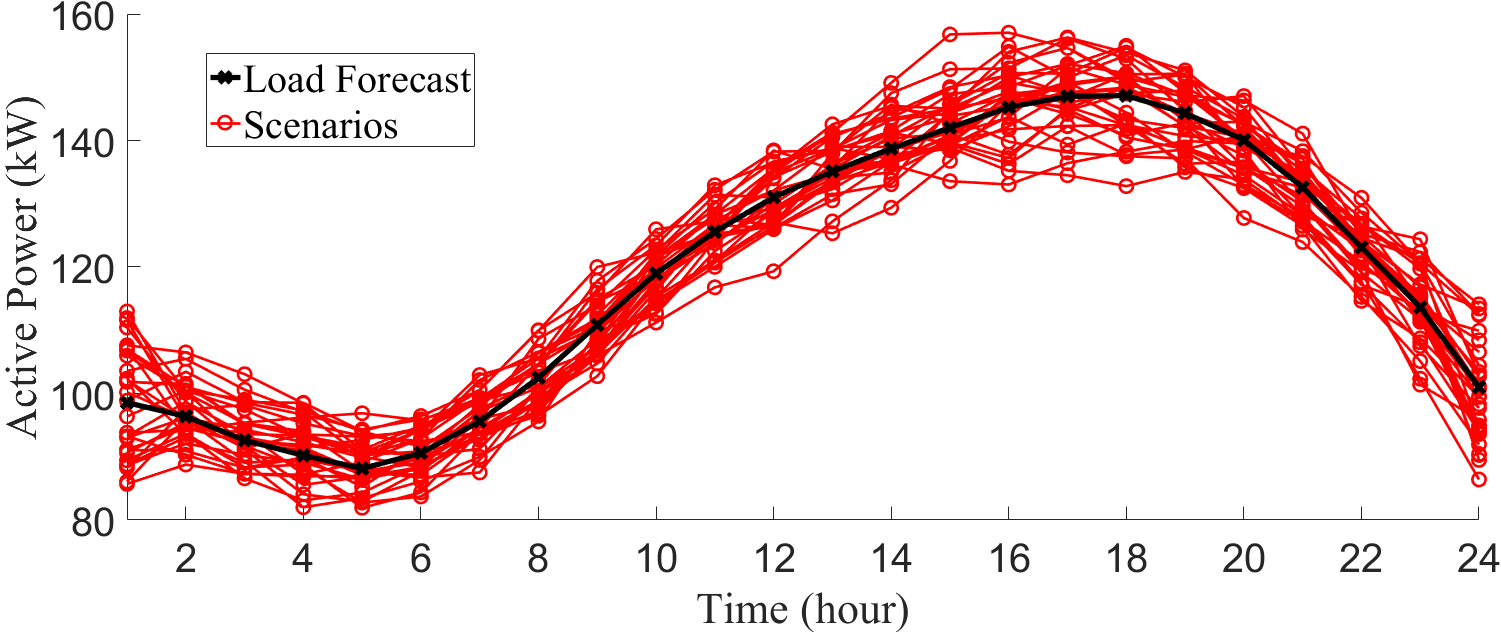}
    \caption{Generated scenarios of active power of a load.}\label{LoadScenarios}
\end{figure}

\vspace{0.2cm}

Each damaged component $m$ is characterized by the repair time $\mathcal{T}_{m,s}$ in scenario $s$. Define $\bm{\mathcal{T}}_s=[\mathcal{T}_{1,s},\mathcal{T}_{2,s},\mathcal{T}_{3,s},..., \mathcal{T}_{D,s}]\in \mathbb{R}^D$ as the vector of real numbers representing the repair time for each damaged component in scenario $s$, where $D$ is the number of damaged components. For $I$ loads and time horizon $T$, let $\bm{e}_s=[e_{1,1,s},e_{1,2,s},...,e_{1,T,s},e_{2,1,s},...,
e_{2,T,s},...,e_{I,1,s},...,e_{I,T,s}]\in \mathbb{R}^{I \cdot T}$ represent the load forecast error in each time period in scenario $s$. By combining $\bm{\mathcal{T}}_s$ and $\bm{e}_s$, the number of random variables is $D+I \cdot T$, and we assume they are mutually independent. Therefore, for $|S|$ scenarios, we can define a matrix $\xi \in \mathbb{R}^{D+I \cdot T\times |S|}$ whose rows consist of random variables and columns consist of scenarios as follows:

{
\small
\[
\bm{\xi} = 
\begin{blockarray}{cccccc}
s=1 & s=2 & s=3 & \dots & s=|\mathcal{S}| \\    \begin{block}{(ccccc)l}
  \mathcal{T}_{1,1} & \mathcal{T}_{1,2} & \mathcal{T}_{1,3} & \dots & \mathcal{T}_{1,|\mathcal{S}|} & v=1 \\      \mathcal{T}_{2,1} & \mathcal{T}_{2,2} & \mathcal{T}_{2,3} & \dots & \mathcal{T}_{2,|\mathcal{S}|} & v=2 \\       \vdots & \vdots & \vdots & \ddots & \vdots & \vdots \\      \mathcal{T}_{D,1} & \mathcal{T}_{D,2} & \mathcal{T}_{D,3} & \dots & \mathcal{T}_{D,|\mathcal{S}|} & v=D \\
e_{1,1,1} & e_{1,1,2} & e_{1,1,3} & \dots & e_{1,1,|\mathcal{S}|} &  v=D+1\\ 
e_{1,2,1} & e_{1,2,2} & e_{1,2,3} & \dots & e_{1,2,|\mathcal{S}|} &  v=D+2\\ 
 \vdots & \vdots & \vdots & \ddots & \vdots & \vdots \\
 e_{1,T,1} & e_{1,T,2} & e_{1,T,3} & \dots & e_{1,T,|\mathcal{S}|} &  v=D+T\\
 e_{2,1,1} & e_{2,1,2} & e_{2,1,3} & \dots & e_{2,1,|\mathcal{S}|} &  v=D+T+1\\ 
e_{2,2,1} & e_{2,2,2} & e_{2,2,3} & \dots & e_{2,2,|\mathcal{S}|} &  v=D+T+2\\ 
 \vdots & \vdots & \vdots & \ddots & \vdots & \vdots \\
 e_{2,T,1} & e_{2,T,2} & e_{2,T,3} & \dots & e_{2,T,|\mathcal{S}|} &  v=D+2T\\
 \vdots & \vdots & \vdots & \ddots & \vdots & \vdots \\
 e_{I,T,1} & e_{I,T,2} & e_{I,T,3} & \dots & e_{I,T,|\mathcal{S}|} &  v=D+I~T\\ \end{block}
\end{blockarray}
 \]
 }where $\xi_{v,s}$ is the realization of random variable $v$ in scenario $s$. According to the Monte Carlo sampling procedure, the probability Pr(s) of each scenario is 1/$|\mathcal{S}|$.
}
\section{Mathematical Formulation}\label{chap:3}

The repair and restoration problem can be divided into two stages. The first stage is to  route the repair crews, which is characterized by depots, repair crews, damaged components and paths between the damaged components. The second stage is distribution system restoration using DGs and reconfiguration. In practice, these two subproblems are interdependent. Therefore, we propose a single MILP formulation that integrates the two problems for joint distribution system repair and restoration, with the objective of maximizing the picked-up loads. {\color{black}The utility solves the optimization problem to obtain the best route for the repair crews. The crews are then dispatched to repair the damaged components. For example, the crews may have to replace a pole or reconnect a wire. This repair process is included in the model through the repair time. Meanwhile, the utility controls the DGs and switches to restore power to the consumers.}


\vspace{-0.5cm}
\subsection{First Stage: Repair Crew Routing}

The routing problem can be defined by a complete graph with nodes and edges $\mathcal{G}(N,E)$. The node set $N$ in the undirected graph contains the depot and damaged components, and the edge set ${E = \left\lbrace(m,n) | m,n \in N; m \neq n \right\rbrace}$ represents the edges connecting each two components. 
Our purpose is to find an optimal route for each crew to reach the damaged components. The value of $x_{m,n,c}$ determines whether the path crew $c$ travels includes the edge $(m,n)$ with $m$ preceding $n$. The routing constraints for the first stage problem are formulated as follows:

\begin{equation}
\mathop \sum \limits_{\forall m \in N}{x_{o_c,m,c}} = 1, \forall c
\label{start_from_depot}
\end{equation}
\begin{equation}
\mathop \sum \limits_{\forall m \in N} {x_{m,d_c,c}} = 1, \forall c
\label{back_to_depot}
\end{equation}
\begin{equation}
\mathop \sum_{\mathclap{\forall n \in N\backslash \left\{ m \right\}}} {x_{m,n,c}} - \mathop \sum_{\mathclap{{\forall n \in N\backslash \left\{ m \right\}}}} {x_{n,m,c}} = 0\;,\;\forall c,\;m \in N\backslash \left\{o_c,d_c\right\}
\label{x=0}
\end{equation}
\begin{equation}
\mathop \sum \limits_{\forall c} \mathop \sum \limits_{\forall m \in N\backslash \left\{n\right\}}x_{m,n,c}=1, \forall n \in N\backslash \left\{o_c,d_c\right\}
\label{x=1}
\end{equation}
{\color{black}
Constraints (\ref{start_from_depot}) and (\ref{back_to_depot}) guarantee that each crew starts and ends its route at the defined start and end locations. For example, if crew 1 is located at the depot, then $x_{o_c,2,1}$=1 means that crew 1 travels from the depot to the damaged component 2. Constraint (\ref{x=0}) is known as the flow conservation constraint; i.e., once a crew repairs the damaged component, the crew moves to the next location. Constraint (\ref{x=1}) ensures that each damaged component is repaired by only one of the crews.
}



\vspace{-0.5cm}
\subsection{Second Stage: Distribution Network Operation}

\subsubsection{Objective}

\begin{equation}
\max \sum \limits_{\forall s } \mathop \sum \limits_{\forall t} \mathop \sum \limits_{\forall i} \textrm{Pr(s)}~ {\omega_i} y_{i,t,s} P^D_{i,t,s}
\label{sec_obj}
\end{equation}
The objective (\ref{sec_obj}) of the second stage is to maximize the expected priority-weighted served loads over the time horizon. In this paper, we consider two load priorities levels: high and low \cite{pr_load}. Note that load priorities can be changed by the utilities as desired. The method in \cite{pr_load} is used to calculate the weights for each load. In the second stage, DGs and line switches are optimally operated in response to the realization of the repair times. Once a damaged line is repaired and energized, it provides a path for the power flow.

\subsubsection{Cold Load Pickup (CLPU)}

\begin{figure}[h!]
\vspace{-0.2cm}
\setlength{\abovecaptionskip}{0pt} 
\setlength{\belowcaptionskip}{0pt} 
  \centering
    \includegraphics[width=0.49\textwidth]{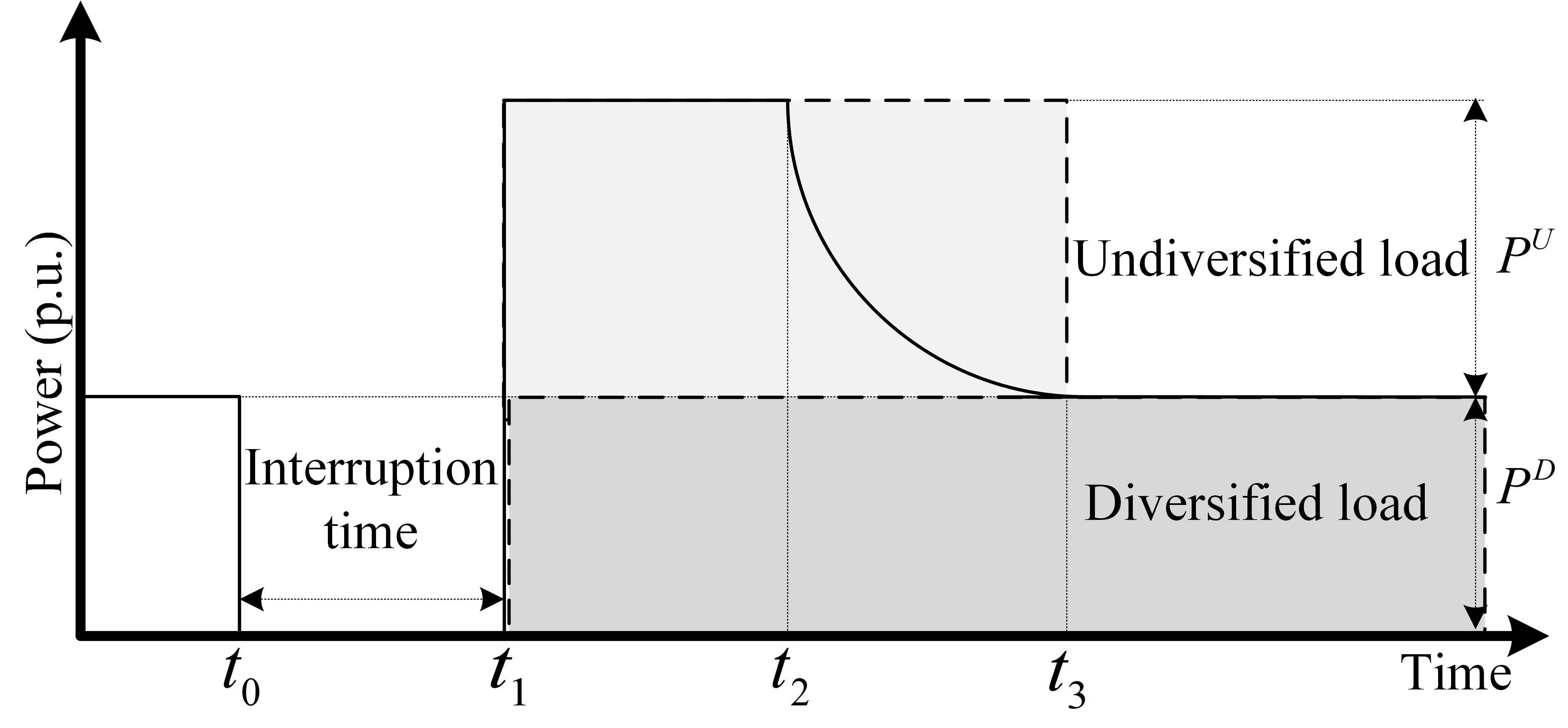}
    \caption{CLPU condition as a delayed exponential model, and the shaded areas represent the two-block model.}\label{Combined_CLPU}
\end{figure}

After an extended period of outage, the effect of cold load pick-up (CLPU) may happen, which is caused by the loss of diversity and simultaneous operation of thermostatically controlled loads. As depicted in Fig. \ref{Combined_CLPU}, the normal steady-state load consumption is defined as the diversified load, and undiversified load is the startup load consumption upon restoration. {\color{black}The time when the load experiences an outage is $t_0$, $t_1$ is the time when the load is restored, and $t_3$ is the time when the load returns to normal condition. The typical behavior of CLPU can be represented using a delayed exponentially decaying function \cite{Liu_2009}, which is shown in Fig. \ref{Combined_CLPU}, where $t_2-t_1$ is the exponential decay delay, and $t_3-t_1$ is the CLPU duration.} This exponential function can be approximated using a linear combination of multiple blocks.

In this paper, we employ two blocks to represent CLPU as suggested in \cite{Liu_2009}. The first block is for the undiversified load $P^U$ and the second for the diversified load $P^D$ (i.e., the steady-state load consumption) as shown in Fig. \ref{Combined_CLPU}. The use of two blocks decreases the computational burden imposed by nonlinear characteristics of CLPU and provides a conservative approach to guarantee the supply-load balance. {\color{black}For a time horizon $T$ and time step $\Delta t$, the CLPU curve is sampled as shown in Fig. \ref{CLPU_3}, where $\lambda$ is the number of time steps required for the load to return to normal condition. The value of $\lambda$ equals the CLPU duration divided by the time step. The CLPU constraint for active power can be formulated as follows:
\begin{equation}
P_{i,t,s}^L=y_{i,t,s}P^D_{i,t,s}+(y_{i,t,s}-y_{i,\rm{max}(t-\lambda,0),s})P^U_{i,t,s},\forall i,t,s
\label{P_clpu}
\end{equation}
where $y_{i,0,s}$ is the initial state of load $i$ immediately after an outage event; i.e., $y_{i,0,s}=1$ and $P_{i,0,s}^L=P_{i,0,s}^D$ if the load is not affected by the outage. If a load goes from a de-energized state to an energized state at time step $t=h$ ($y_{i,h-1,s}=0$ and $y_{i,h,s}=1$), it will return to normal condition at time step $h+\lambda$, as $y_{i,h,s}-y_{i,max(h+\lambda-\lambda,0),s}=0$. Before time step $h+\lambda$, $P_{i,t,s}^U$ is added to $P_{i,t,s}^D$ to represent the undiversified load. The function max($t-\lambda,0$), is used to avoid negative values. We assume that the duration of the CLPU decaying process is one hour in the simulation \cite{Liu_2009}. Moreover, the study in \cite{Nagpal_2014} showed that the total load at pick-up time can be up to 200\% of the steady state value, thus, $P_{i,t,s}^U$ is set to be equal to $P_{i,t,s}^D$. Similarly, the CLPU constraint for reactive power can be formulated as follows:
\begin{multline}
Q_{i,t,s}^L=y_{i,t,s}Q^D_{i,t,s}+\\(y_{i,t,s}-y_{i,\rm{max}(t-\lambda,0),s})Q^U_{i,t,s},\forall i,t,s
\label{Q_clpu}
\end{multline}

\begin{figure}[h!]
\vspace{-0.2cm}
\setlength{\abovecaptionskip}{0pt} 
\setlength{\belowcaptionskip}{0pt} 
  \centering
    \includegraphics[width=0.49\textwidth]{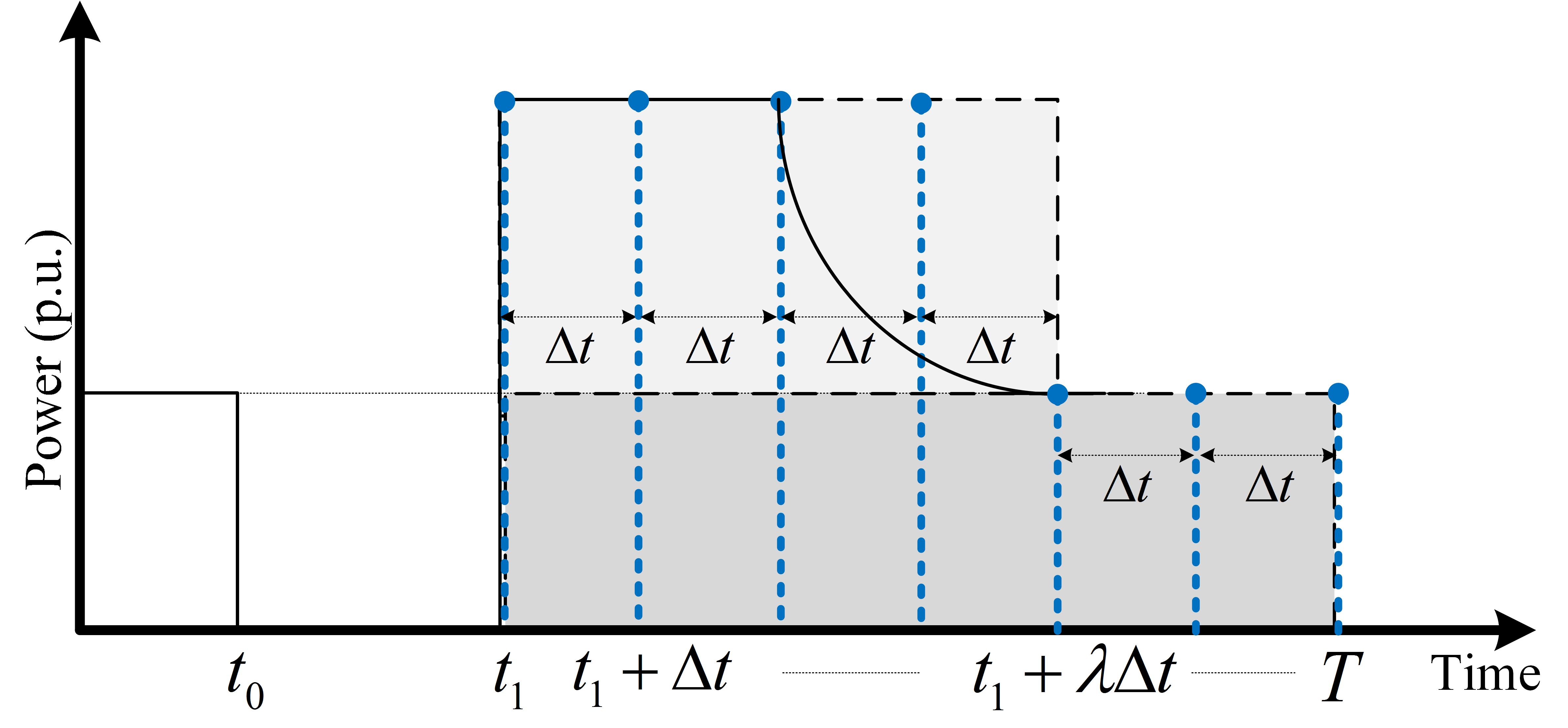}
    \caption{Two-blocks CLPU condition as a delayed exponential model, with time step $\Delta t$.}\label{CLPU_3}
\end{figure}
}
\subsubsection{Distribution Network Optimal Power Flow}

The power flow model mostly used in transmission network restoration is the linear DC optimal power flow model which neglects reactive power and voltage levels. AC optimal power flow, on the other hand, is nonlinear and will greatly increase the computational burden of the problem. Therefore, linearized Distflow equations are used to calculate the power flow and the voltages at each node. Linearized Distflow equations have been used and verified in the literature \cite{distflow2,Arif2017,Ma2017,Wang2015}. The equations are formulated as follows:

\begin{equation}
\mathop \sum_{\mathclap{\forall k \in K\left( {.,i} \right)}} {P^B_{k,t,s}} + P^G_{i,t,s}= \mathop \sum_{\mathclap{\forall k \in K\left( {i,.} \right)}} {P^B_{k,t,s}} + P_{i,t,s}^L,\forall i,t,s
\label{power balance}
\end{equation}
\begin{equation}
\mathop \sum_{\mathclap{\forall k \in K\left( {.,i} \right)}} {Q^B_{k,t,s}} + Q_{i,t,s}^{G} = \mathop \sum_{\mathclap{\forall k \in K\left( {i,.} \right)}} {Q^B_{k,t,s}} + Q_{i,t,s}^L,\forall i,t,s
\label{reactive balance}
\end{equation}
\begin{equation}
{V_{j,t,s}} - {V_{i,t,s}} + \frac{{{R_k}{P^B_{k,t,s}} + {X_k}{Q^B_{k,t,s}}}}{{{V_1}}}\;\le (1-u_{k,t,s}) M,\forall k,t,s
\label{voltage1}
\end{equation}
\begin{equation}
(u_{k,t,s}-1) M \;\le\; {V_{j,t,s}} - {V_{i,t,s}} + \frac{{{R_k}{P^B_{k,t,s}} + {X_k}{Q^B_{k,t,s}}}}{{{V_1}}},\forall k,t,s
\label{voltage2}
\end{equation}
\begin{equation}
1-\epsilon  \le {V_{i,t,s}} \le 1+\epsilon \;,\;\forall i,t,s
\label{V limits}
\end{equation}
Constraints (\ref{power balance}) and (\ref{reactive balance}) represent the active and reactive power balance constraints, respectively. The voltage at each bus is expressed in constraints (\ref{voltage1}) and (\ref{voltage2}), where $V_1$ is the reference voltage. A disjunctive method is used to ensure that the voltage levels of two disconnected buses are decoupled. The values used for $M$ are explained in Section \ref{M_section}. Constraint (\ref{V limits}) defines the allowable range of voltage deviations, where $\epsilon$ is set to be 5$\%$.

We consider dispatchable DGs for supplying loads in the distribution network, and automatic switches to reconfigure the network. The automatic switches are controlled by $u_{k,t,s}, k \in \textcolor{black}{\Omega_{SW}}$. The following constraints define the capacity of the DGs, line flow limits, and switching status of the lines:

\begin{equation}
0 \le P_{i,t,s}^{G} \le P_i^{G_{max}}\;,\;\forall i,t,s
\label{DG limits}
\end{equation}
\begin{equation}
0 \le Q_{i,t,s}^{G} \le Q_i^{G_{max}}\;,\;\forall i,t,s
\label{QDG limits}
\end{equation}
\begin{equation}
- u_{k,t,s}P_k^{B_{max}} \le {P^B_{k,t,s}} \le u_{k,t,s}P_k^{B_{max}}\;,\;\forall k,t,s
\label{PL limits}
\end{equation}
\begin{equation}
- u_{k,t,s}Q_k^{B_{max}} \le {Q^B_{k,t,s}} \le u_{k,t,s}Q_k^{B_{max}}\;,\;\forall k,t,s
\label{QL limits}
\end{equation}
\begin{equation}
u_{k,t,s} = 1, \forall k \not \in \{\Omega_{SW} \cup N\backslash\{0\}\},s
\label{uk=1}
\end{equation}

Constraints (\ref{DG limits}) and (\ref{QDG limits}), respectively, define the real and reactive output limits for DGs. Constraints (\ref{PL limits}) and (\ref{QL limits}) set the limits of the line flows and indicate that the power flow through a damaged line equals zero, which is achieved by multiplying the line limits by $u_{k,t,s}$. Constraint (\ref{uk=1}) maintains the switching status of a line $u_{k,t,s}$ to be 1 when there is no damage and/or no switch.

Once a load is served, it should remain energized, as enforced by the following constraint:

\begin{equation}
{y_{i,t+1,s}} \ge {y_{i,t,s}}\;,\;\forall i,t,s
\label{units on}
\end{equation}

\subsubsection{Radiality Constraints}

The distribution network is reconfigured dynamically using switches to change the topology of the network. Radiality constraints are introduced to maintain radial configuration. The method used in \cite{Borghetti2012} is employed in this paper. Radiality is enforced by introducing constraints for ensuring that at least one of the lines of each possible loop in the network is open. A depth-first search method \cite{Borghetti2012} is used to identify the possible loops in the network and the lines associated with them. The following constraint can then be used to ensure radial configuration:

\begin{equation}
\sum_{k \in \Omega_{K(l)}} u_{k,t,s} \le |\Omega_{K(l)}|-1, \forall l,t,s
\label{radial_con}
\end{equation}
where $|\Omega_{K(l)}|$ is the number of lines in loop $l$. Constraint (\ref{radial_con}) guarantees that at least one line is disconnected in each loop. {\color{black}Alternatively, the radiality constraints can be represented by (\ref{beta_limits})-(\ref{oneparrent}) based on the spanning tree approach \cite{recon}.

\begin{equation}
0 \le \beta^s_{i,j,t} \le 1, \forall i,j \in \Omega_B,t,s 
\label{beta_limits}
\end{equation}
\begin{equation}
\beta^s_{i,j,t}+\beta^s_{j,i,t}=u_{k,t,s},\; \forall k,t,s
\label{radialconnection}
\end{equation}
\begin{equation}
\beta^s_{i,j,t}=0,\; \forall i \in \Omega_B, j \in \Omega_{SB},t,s
\label{substation}
\end{equation}
\begin{equation}
\mathop \sum \limits_{\forall i \in \Omega_B} \beta^s_{i,j,t} \le 1,\;\forall j \in \Omega_B,t,s
\label{oneparrent}
\end{equation}
Two variables $\beta_{i,j,t}$ and $\beta_{j,i,t}$ are defined to model the spanning tree. For a radial network, each bus cannot be connected to more than one parent bus and the number of lines equals the number of buses other than the root bus. Constraint (\ref{radialconnection}) relates the connection status of the line and the spanning tree variables $\beta_{i,j,t}$ and $\beta_{j,i,t}$. If the distribution line is connected, then either $\beta_{i,j,t}$ or $\beta_{j,i,t}$ must equal one. Constraint (\ref{substation}) designates substations as and indicates that they do not have parent buses. Constraint (\ref{oneparrent}) requires that every bus has no more than one parent bus. The spanning tree constraints guarantee that the number of buses in a spanning tree, other than the root, equals the number of lines \cite{recon}. In this paper, we use constraint (\ref{radial_con}) to ensure the radiality as the spanning tree constraints in (\ref{beta_limits})-(\ref{oneparrent}) will add $|\Omega_{B}|\times|\Omega_{B}|\times|T|\times|S|$ variables.
}
\subsubsection{Restoration Time}

The arrival time and consequently the time when each component is repaired must be calculated to connect the routing and power operation problems. Once a crew arrives at a damaged component $m$ at time $\alpha_{m,c}$, they spend a time $\mathcal{T}_{m,s}$ to repair the damaged component, and then take time $T^R_{m,n,c}$ to arrive at the next damaged component $n$. Therefore, $\alpha_{m,c,s} + {\mathcal{T}_{m,s}} + T^R_{m,n} = \alpha_{n,c,s}$ if crew $c$ travels the path $m$ to $n$. The travel time between the damaged components and depot can be obtained through a geographic information system (GIS). The arrival time constraints are formulated as follows:
\vspace{-0.1cm}
\begin{equation}
\begin{gathered}
\alpha_{m,c,s} + {\mathcal{T}_{m,s}} + T^R_{m,n} - \left( {1 - {x_{m,n,c}}} \right)M \le \alpha_{n,c,s}\\\;\forall m \in N \backslash \{d_c\},n \in N\backslash \left\{o_c,m\right\},c,s
\end {gathered}
\label{AT1}
\end{equation}
\begin{equation}
\begin{gathered}
 \alpha_{n,c,s} \le \alpha_{m,c,s} + {\mathcal{T}_{m,s}} + T^R_{m,n}+\left( {1 - {x_{m,n,c}}} \right)M\\\;\forall m \in N \backslash \{d_c\},n \in N\backslash \left\{o_c,m\right\},c,s
\end {gathered}
\label{AT1b}
\end{equation}

Disjunctive constraints are used to decouple the times to arrive at components $m$ and $n$ if the crew does not travel from $m$ to $n$. In order to determine when will the damaged component be restored and can be operated again, we enforce the following constraints: 

\begin{equation}
0 \le f_{m,t,s} \le 1, \forall m \in N\backslash \left\{ {o_c,d_c} \right\},t,s
\end{equation}
\begin{equation}
\mathop \sum \limits_{\forall t} {f_{m,t,s}} = 1\;,\;\forall m \in N\backslash \left\{ {o_c,d_c} \right\},s
\label{f=1}
\end{equation}
For example, if component $m$ is repaired at $t=3$, then $f_m=\{0,0,1,0,...,0\}$. The restoration time for component $m$ can be found by $\mathop \sum_{\forall t}t\;{f_{m,t}}$.

The restoration time depends on the arrival time and the repair time, where the relationship is modeled using the following equations:
\vspace{-0.2cm}
\begin{equation}
\begin{gathered}
\mathop \sum \limits_{\forall t} t{f_{m,t,s}} \ge \mathop \sum \limits_{\forall {c}} \left( \alpha_{m,c,s} + {\mathcal{T}_{m,s}}\mathop \sum \limits_{\forall n \in N } x_{m,n,c} \right)\\ \forall m \in N\backslash \left\{ {o_c,d_c} \right\},s
\end{gathered}
\label{tf1}
\end{equation}
\begin{equation}
\begin{gathered}
\mathop \sum \limits_{\forall t} t{f_{m,t,s}} \le \mathop \sum \limits_{\forall {c}} \left( \alpha_{m,c,s} + {\mathcal{T}_{m,s}}\mathop \sum \limits_{\forall n \in N } x_{m,n,c} \right)\\+1-\epsilon, \forall m \in N\backslash \left\{ {o_c,d_c} \right\},s
\end{gathered}
\label{tf1b}
\end{equation}
\begin{equation}
0 \le \alpha_{m,c,s} \le M \mathop \sum \limits_{n \in N}x_{m,n,c},\;\forall m \in N\backslash \left\{ {o_c,d_c} \right\},c,s
\label{AT2}
\end{equation}

\begin{figure}[h!]
\vspace{-0.2cm}
\setlength{\abovecaptionskip}{0pt} 
\setlength{\belowcaptionskip}{0pt} 
  \centering
    \includegraphics[width=0.49\textwidth]{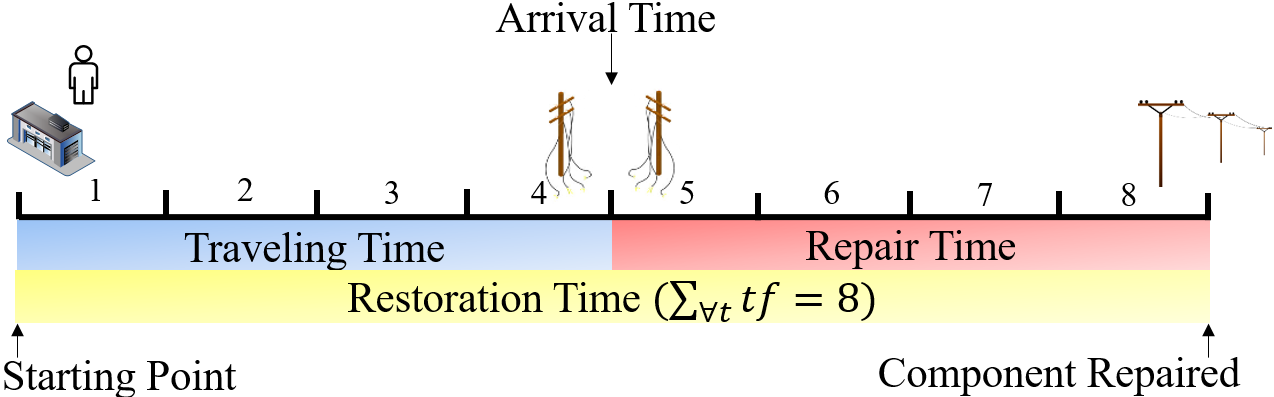}
    \caption{Time sequence of the repair process.}\label{time_label}
\end{figure}

Constraints (\ref{tf1}) and (\ref{tf1b}) determine the time when a damaged component is repaired by adding its repair time to the arrival time. The two equations are used to define $\ceil{tf_{m,t}}$, since the time horizon has integer values. If the damaged component is not repaired by a crew $c$, then the arrival time and repair time for this crew should not affect constraints (\ref{tf1}) and (\ref{tf1b}), which is realized by using constraint (\ref{AT2}) to set $\alpha_{m,c}=0$. Fig. \ref{time_label} demonstrates the time sequence of the repair process and how to find the restoration time. Starting from the depot, if both travel time and repair time are 4 hours, the restoration time is $\mathop \sum \limits_{\forall t}t\;{f_{m,t}}=8$.

The routing and power operation problems are connected with the following constraint:

\begin{equation}
{u_{m,t,s}} = \mathop \sum \limits_{\bar{t}  = 1}^{t} {f_{m,\bar{t},s}}\;,\;\forall m \in N\backslash \left\{ {o_c,d_c} \right\},t,s
\label{z=f}
\end{equation}
Constraint (\ref{z=f}) indicates that the restored component becomes available after it is repaired, and remains available in all subsequent time periods. We assume that the repair time includes the time it takes to re-energize the component; therefore, if component $m$ is repaired at $t$ = 4, it can be used at $t$ = 4 and thereafter. For example, if $t=[1,2,...,6]$ and $f_m=[0,0,0,1,0,0]$ then $u_{m,t}=[0,0,0,1,1,1]$.

{\color{black}
\subsubsection{Big M} \label{M_section}

The value used for $M$ depends on the constraint. An inappropriately large $M$ may increase the computation time, and a small value may introduce infeasibility. In constraint (\ref{voltage1}) and (\ref{voltage2}), the maximum and minimum values for the voltage are 1.05 and 0.95 per unit. Hence, the largest possible difference between any two voltages ($V_{j,t,s}-V_{i,t,s}$) is 0.1 per unit. Also, the maximum drop in voltage ($R_k P_{k,t,s}^B+X_k Q_{k,t,s}^B)/V_1$ is 0.1 per unit. Accordingly, the minimum value of $M$ in (\ref{voltage1}) and (\ref{voltage2}) is 0.2 per unit.

In the routing constraints, the crews must arrive at the damaged components before starting the repairs. For example, if the time horizon  is $T=10$, and the repair time for some damaged component $m$ is  $T_{m,s}=1$, then the crew should arrive at $\alpha_{m,c,s}=9$ at the latest in order to repair the component. Note that the time horizon should be chosen such that all damaged components can be repaired in the optimization problem. Therefore, the minimum value of $M$ in (\ref{AT2}) equals the time horizon minus the minimum repair time. The minimum repair time is used to obtain the largest difference between $T$ and the repair times of the components. Denote the value of $M$ in (\ref{AT2}) as $M_{27}$. For (\ref{AT1}) and (\ref{AT1b}), the value of $M$ should be larger than the time horizon $T$. In a worst-case scenario, the arrival time of crew $c$ at damaged component $m$ is $\alpha_{m,c,s}=M_{27}$, and the crew does not repair damaged component $n$, as per equation (\ref{AT2}), $\alpha_{n,c,s}=0$. Consequently, (\ref{AT1}) and (\ref{AT1b}) are translated to $-M \le 0-M_{27}-\mathcal{T}_{m,s}-T_{m,n}^R \le M$. Hence, the minimum value of $M$ in (\ref{AT1}) and (\ref{AT1b}) equals $M_{27}$ plus the maximum repair and travel times.
}

\subsection{Two-Stage Stochastic Program}

In this paper, we formulate the stochastic DSRRP as a two-stage stochastic program. In the first stage, the crews are dispatched to the damaged components. Therefore, the first-stage variable is $x_{m,n,c}$. After realization of the repair times and loads, the distribution network is operated in the second stage. The second-stage variables are defined in vector $\bm{\gamma_s}$, which includes $(\alpha,f,P^B,P^G,P^L,Q^B,Q^G,Q^L,u,V,y,\beta)$. The \textit{extensive form} (EF) of the two-stage stochastic DSRRP is formulated as follows:

\vspace{-0.2cm}

\begin{equation}
\zeta  \textcolor{black}{(\textrm{weighted kWh})}= \max_{\bm{x},\bm{\gamma}}\sum \limits_{\forall s } \mathop \sum \limits_{\forall t} \mathop \sum \limits_{\forall i} \textrm{Pr(s)} {\omega_i} y_{i,t,s}P^D_{i,t,s}
\label{ef_obj}
\end{equation}
\begin{center}
s.t. (\ref{start_from_depot})-(\ref{x=1}), (\ref{P_clpu})-(\ref{z=f})
\end{center}
\begin{equation}
u,x,y \in \left\{ {0,1} \right\} 
\label{binary variables}
\end{equation}

\vspace{-0.2cm}

\section{Solution Algorithm}\label{chap:4}

In this section, we decompose S-DSRRP and present the algorithm for solving the decomposed problem.

\subsection{Progressive Hedging} 

{\color{black}Watson and Woodruff adapted the PH algorithm \cite{Rockafellar1991} to approximately solve stochastic mixed-integer problems. The PH algorithm decomposes the extensive form into subproblems, by relaxing the non-anticipativity of the first-stage variables. Hence, for $|\mathcal{S}|$ scenarios, the stochastic program is decomposed into $|\mathcal{S}|$ subproblems. PH can solve the subproblems in parallel to reduce the computational burden for large-scale instances. The authors of \cite{Ryan2013} effectively implemented PH for solving the stochastic unit commitment problem. A full description of the PH algorithm can be found in \cite{Rockafellar1991}.

To demonstrate the PH algorithm, we first define a compact form for the general two-stage stochastic program as follows:

\begin{equation}
\zeta = \min_{\bm{\delta,\gamma_s}} ~~~\bm{a}^T \bm{\delta} + \sum \limits_{\forall s }\textrm{Pr(s)} ~ \bm{b_s}^T \bm{\gamma_s}
\end{equation}
\begin{equation}
\textrm{s.t.} ~~(\bm{\delta},\bm{\gamma_s})\in \mathcal{Q}_s, \forall s
\end{equation}
where \bm{$a$} and \bm{$b_s$} are vectors containing the coefficients associated with the first-stage (\bm{$\delta$}) and second-stage (\bm{$\gamma_s$}) variables in the objective, respectively. The restriction $(\bm{\delta,\gamma_s})\in \mathcal{Q}_s$ represents the subproblem constraints that ensures a feasible solution. The PH algorithm is described in Algorithm 1, using a penalty factor $\rho$ and a termination threshold $\varepsilon$.

 \begin{algorithm}[h!]
 	\caption{The Two-Stage PH Algorithm}
 	\small
 	\begin{algorithmic}[1]
 		\STATE Let $\tau:=0$
 		\STATE For all  $s \in \mathcal{S}$, compute:\\ $\bm{\delta_s}^{(\tau)} := \argmin_{\bm{\delta}} \Big{\{} \bm{a}^T \bm{\delta} + \bm{b_s}^T \bm{\gamma_s} : (\bm{\delta},\bm{\gamma_s})\in \mathcal{Q}_s \Big{\}}$
 	\STATE $\bm{\bar{\delta}}^{(\tau)}:=\sum_{s \in \mathcal{S}}\textrm{Pr(s)}\bm{\delta_s}^{(\tau)}$
 	\STATE $\bm{\eta_s}^{(\tau)}:=\rho(\bm{\delta_s}^{(\tau)}-\bm{\bar{\delta}}^{(\tau)})$
 	\STATE $\tau:=\tau+1$
 	\STATE For all $s \in \mathcal{S}$ compute:\\
 	$\bm{\delta_s}^{(\tau)}:=\argmin_{\bm{\delta}} \Big{\{} \bm{a}^T \bm{\delta} + \bm{b_s}^T \bm{\gamma_s}+ \bm{\eta_s}^{(\tau-1)}\bm{\delta} + \frac{\rho}{2}||{\bm{\delta}-\bm{\bar{\delta}}^{(\tau-1)}}||^2: (\bm{\delta},\bm{\gamma_s})\in \mathcal{Q}_s \Big{\}}$
 	\STATE $\bm{\bar{\delta}}^{(\tau)}:=\sum_{s \in \mathcal{S}}\textrm{Pr(s)}\bm{\delta_s}^{(\tau)}$
 	\STATE $\bm{\eta_s}^{(\tau)}:=\bm{\eta_s}^{(\tau-1)}+\rho(\bm{\delta_s}^{(\tau)}-\bm{\bar{\delta}}^{(\tau)})$
 	\STATE $\mu^{(\tau)}:=\sum_{s \in \mathcal{S}}\textrm{Pr(s)}||\bm{\delta_s}^{(\tau)}-\bm{\bar{\delta}}^{(\tau)}||$
 	\STATE If $\mu^{(\tau)}<\textcolor{black}{\varepsilon}$, then go to \textbf{Step 5}. Otherwise, terminate
 	\end{algorithmic}
 \end{algorithm}
 }
 
The PH algorithm starts by solving the subproblems with individual scenarios in Step 2. Notice that for an individual scenario, the two-stage model boils down to a single-level problem. Step 3 aggregates the solutions to obtain the expected value $\bm{\bar{\delta}}$. The multiplier $\eta_s$ is updated in Step 4. The first four steps represent the initialization phase. In Step 6, the subproblems are augmented with a linear term proportional to the multiplier $\bm{\eta}^{(\tau-1)}_s$ and a squared two norm term penalizing the difference of $\bm{\delta}$ from $\bm{\bar{\delta}}^{(\tau-1)}$, where $\tau$ is the iteration number. Steps 7-8 repeat Steps 3-4. The program terminates once $\sum_{s \in \mathcal{S}}\textrm{Pr(s)}||\bm{\delta_s}^{(\tau)}-\bm{\bar{\delta}}^{(\tau)}|| < \varepsilon$; i.e., all first-stage decisions $\bm{\delta_s}$ converge to a common $\bm{\bar{\delta}}$. The termination threshold $\varepsilon$ is set to be 0.01 in this paper.

 \subsection{Decomposed S-DSRRP}

The proposed algorithm iteratively selects a group of damaged components and dispatches the crews until all damaged components are repaired. The S-DSRRP is decomposed into two subproblems.

\subsubsection{Subproblem I} 
The first subproblem determines $\mathcal{C}$ critical damaged components to repair. This problem is formulated as a two-stage SMIP. In the first stage, the critical damaged components are determined, and the distribution network is operated in the second stage. The first subproblem is formulated as follows:
\vspace{-0.2cm}
\begin{equation}
  \bm{z}^* := \argmax_{\bm{z},\bm{\bar{\gamma_s}}}~ \sum \limits_{\forall s } \mathop \sum \limits_{\forall t} \mathop \sum \limits_{\forall i} \textrm{Pr(s)}~ {\omega_i}~  y_{i,t,s} P^D_{i,t,s}
  \label{Divide_Obj}
  \end{equation}
\begin{center}
s.t. (\ref{P_clpu})-(\ref{radial_con})
\end{center}
  \begin{equation}
\mathop \sum \limits_{\forall m\in N\backslash \{0\}} z_m \le \mathcal{C}
\label{zlenc}
\end{equation}
\begin{equation}
u_{m,t,s} \le z_m, \forall m,t,s
\label{ulez}
\end{equation}
\begin{equation}
\sum_{t=1}^{\mathcal{T}_{m,s}} u_{m,t,s} = 0, \forall m,s
\label{u==0}
\end{equation}where \bm{$\bar{\gamma_s}$} includes $(P^B,P^G,P^L,Q^B,Q^G,Q^L,u,V,y,\beta)$. Define binary variable $z_m$ to equal 1 if damaged component $m$ is a critical damaged component to repair. The goal of this subproblem is to find a number of damaged components that, if repaired, will maximize the served load. In order to obtain a manageable problem for the second subproblem, we set the number of selected (critical) damaged components to be equal to the number of crews; i.e., $\mathcal{C}$. In this subproblem, all routing constraints are neglected, and we assume that the crews instantaneously begin repairing the selected damaged components. The objective of Subproblem I (\ref{Divide_Obj}) is to maximize the served loads, while considering distribution network operation constraints. Constraint (\ref{zlenc}) limits the number of damages to be repaired. 
If $z_m$ equals 0, then $u_{m,t,s}$ must be 0, which is enforced by (\ref{ulez}). Constraint (\ref{u==0}) sets $u_{m,t,s}$ to be 0 until time $\mathcal{T}_{m,s}$ has passed. After determining the critical components, we proceed to the second subproblem.

\subsubsection{Subproblem II}

The second subproblem is formulated similarly to (\ref{ef_obj}). The crews are dispatched to the damaged components obtained from Subproblem I in the first stage, and the distribution network is operated in the second stage. Each cycle of Subproblem I and Subproblem II is defined as a dispatch cycle. The dispatch cycle is denoted by $r$. Define the subset of critical damaged components and starting point as $N'(r)$. Note that the starting point after the first dispatch cycle is the current location of the crew instead of the depot. Subproblem II solves the two-stage S-DSRRP for $N'(r)$, which is formulated as follows:

\begin{center}
$\zeta$ = $\max_{\bm{x},\bm{\gamma_s}}$ (\ref{ef_obj})
\end{center}
\begin{center}
s.t. (\ref{start_from_depot})-(\ref{x=1}), (\ref{P_clpu})-(\ref{radial_con}), (\ref{AT1})-(\ref{z=f})
\end{center}
\begin{equation}
u_{m,t,s}=0, \forall t,s,m \in N\backslash N'(r)
\label{u=0_DC}
\end{equation}

Constraint (\ref{u=0_DC}) states that if component $m$ is damaged and is not being repaired, then $u_{m,t,s}$ equals 0. The two subproblems are repeated until all damaged components are repaired.

Algorithm 2 presents the pseudo-code for the D-PH algorithm. The number of dispatch cycles is equal to the number of damaged components divided by the number of crews; i.e., $\floor{|N\backslash \{\textrm{depot}\}|/\mathcal{C}}$. If there are 11 damages and 3 crews, then the number of dispatch cycles will be 3, and the remaining damaged components are considered in Steps 11-12. The algorithm starts by solving Subproblem I in Step 2 using PH. After obtaining $z^*$ in dispatch cycle $r$, the subset of critical damaged components, $N'(r)$, is defined in Step 3. If $N'(r)$ is null, then all loads can be served without repairing any damaged components. Therefore, the loop ends and the routing problem is solved for $N$ in Step 12. Subproblem II is solved next using PH in Step 7 to route the crews and operate the distribution network. We then update $o_c$ in Step 8 by using the results obtained from the Subproblem II. The end point for the crews is set to be the depot, but the variable $x_{m,d_c,c}$ is used only to determine the starting locations for the next dispatch cycle. The crews return to the depot after all repair tasks are finished in the final dispatch cycle. The set of damaged components is updated in Step 9 by removing the repaired lines. Step 11 checks whether there are any remaining damaged components, and then solves Subproblem II to finish the repairs.

{\color{black}
\begin{algorithm}[h!]
 \caption{D-PH algorithm for solving S-DSRRP}
 \small
 \begin{algorithmic}[1]
 \renewcommand{\algorithmicrequire}{\textbf{Input:}}
 \renewcommand{\algorithmicensure}{\textbf{Output:}}
 \REQUIRE $ \mathcal{C},P_{i,t,s}^D,Q_{i,t,s}^D,\mathcal{T}_{m,s},R_k,X_k,T^R_{m,n},w_i,N$
 \ENSURE  $ \alpha_{m,c,s},P_{i,t,s}^G,Q_{i,t,s}^G,u_{k,t,s},x_{m,n,c},y_{i,t,s}$
  \FOR {$r = 1$ to $\floor{|N\backslash \{\textrm{depot}\}|/\mathcal{C}}$}
  \STATE Solve using \textbf{PH} \COMMENT{Subproblem I}
\\$\bm{z}^* :=\argmax_{\bm{z},\bm{\bar{\gamma_s}}}\{(36):s.t.~(7)\mbox{-}(20),(38)\mbox{-}(40) \}$
  \STATE $N'(r)=\{m|z^*_m = 1, \forall m \in N$\}
  \IF {$N'(r)$ is null}
  \STATE \textcolor{black}{\textbf{break}} \COMMENT{All loads can be served}
  \ENDIF
  \STATE Solve using \textbf{PH} \COMMENT{Subproblem II}\\
$ \zeta :=\max_{\bm{x},\bm{\gamma_s}}\{(33):s.t.(2)\mbox{-}(5),(7)\mbox{-}(20),(25)\mbox{-}(32),(41) \}$
\STATE For each crew, update the starting location:\\ $o_c=\{m|x^*_{m,d_c,c} = 1, \forall m \in N$\}
\STATE $N = N \backslash N'(r)$ \COMMENT{update damaged components}
  \ENDFOR
  \IF {$N$ is not null}
  \STATE Repeat \textbf{Step 7} \COMMENT{route the repair crews to the remaining damaged components}
  \ENDIF
 \end{algorithmic} 
 \end{algorithm}
 }
\vspace{-0.4cm}

\section{Simulation and Results}\label{chap:5}


Modified IEEE 34- and 8500-bus distribution feeders are used as test cases for the repair and restoration problem. Detailed information on the networks can be found in \cite{34ieee} and \cite{8500ieee}, respectively. The stochastic models and algorithms are implemented using the PySP package in Pyomo \cite{Hart2012}. IBM's CPLEX 12.6 mixed-integer solver is used to solve all subproblems. 
The experiments were performed on Iowa State University's Condo cluster, whose individual blades consist of two 2.6 GHz 8-Core Intel E5-2640 v3 processors and 128GB of RAM. The scenario subproblems are solved in parallel by using the Python Remote Objects library. 
To ensure a fast response for the outage, and the convergence of the algorithm, we impose a 30-minute time limit on each subproblem; i.e., a one-hour time limit \cite{last_mile} for each dispatch cycle.

\subsection{Case I: IEEE 34-bus distribution feeder}

The IEEE 34-bus feeder is modified by adding three dispatchable backup DGs installed at randomly selected locations, and two-line switches. High-priority loads are chosen arbitrarily. The capacity of the DGs is 150 kW. The travel time between damaged components ranges from 15 to 30 minutes, and the time step used in the simulation is one hour. We assume three crews, one depot, and seven damaged lines. The outage is assumed to have occurred at 12 AM. 
The Monte Carlo sampling technique is used to generate 1000 random scenarios with equal probability, and the simultaneous backward scenario reduction algorithm \cite{Dupacova2003} is applied to reduce the number of scenarios to 30. The General Algebraic Modeling System (GAMS) provides a toolkit named SCENRED2 for implementing the scenario reduction algorithm \cite{GAMS}. For the repair time, a lognormal distribution is used with parameters $\mu=-0.3072$ and $\sigma=1.8404$ \cite{lognormal2}, and unrealistic values (e.g., 0.01 hours) are truncated. On the other hand, the load forecast error is generated using a truncated normal distribution with limits $\pm$15\% \cite{Lu2013}. Samples of the 30 generated scenarios are shown in Table \ref{RTscenarios} for the repair time.

\vspace{-0.2cm}

\begin{table}[htbp]
\small
  \centering
  \caption{Samples of the repair times (in hours) for the 30 generated scenarios using the lognormal distribution}
  \vspace{-0.25cm}
\begin{tabular}{cccccc}
    \hline
\rule{0pt}{2ex}     Damage & Scenario 1 & Scenario 2 & Scenario 3 & ....  & Scenario 30\\
    \hline
\rule{0pt}{2ex}    Line ~5-6~~~ & 2.71  & 3.61  & 1.97  & ....  & 3.11\\
    Line ~7-8~~ & 4.01  & 2.36  & 3.85  & ....  & 5.11 \\
    Line ~9-10~ & 4.03   & 3.21  & 1.06  & ....  & 4.62 \\
    Line 12-13 & 2.18  & 1.87  & 2.88  & ....  & 3.45 \\
    Line 31-32 & 1.14  & 1.83  & 3.07  & ....  & 6.95 \\
    Line 17-18 & 2.87  & 3.93  & 3.09  & ....  & 8.21 \\
    Line ~4-20 & 1.68  & 1.84  & 4.69  & ....  & 2.46 \\
    \hline
    \end{tabular}%
  \label{RTscenarios}%
\end{table}%
\vspace{0.2cm}

\begin{figure}[h!]
\vspace{-0.2cm}
\setlength{\abovecaptionskip}{0pt} 
\setlength{\belowcaptionskip}{0pt} 
  \centering
    \includegraphics[width=0.49\textwidth]{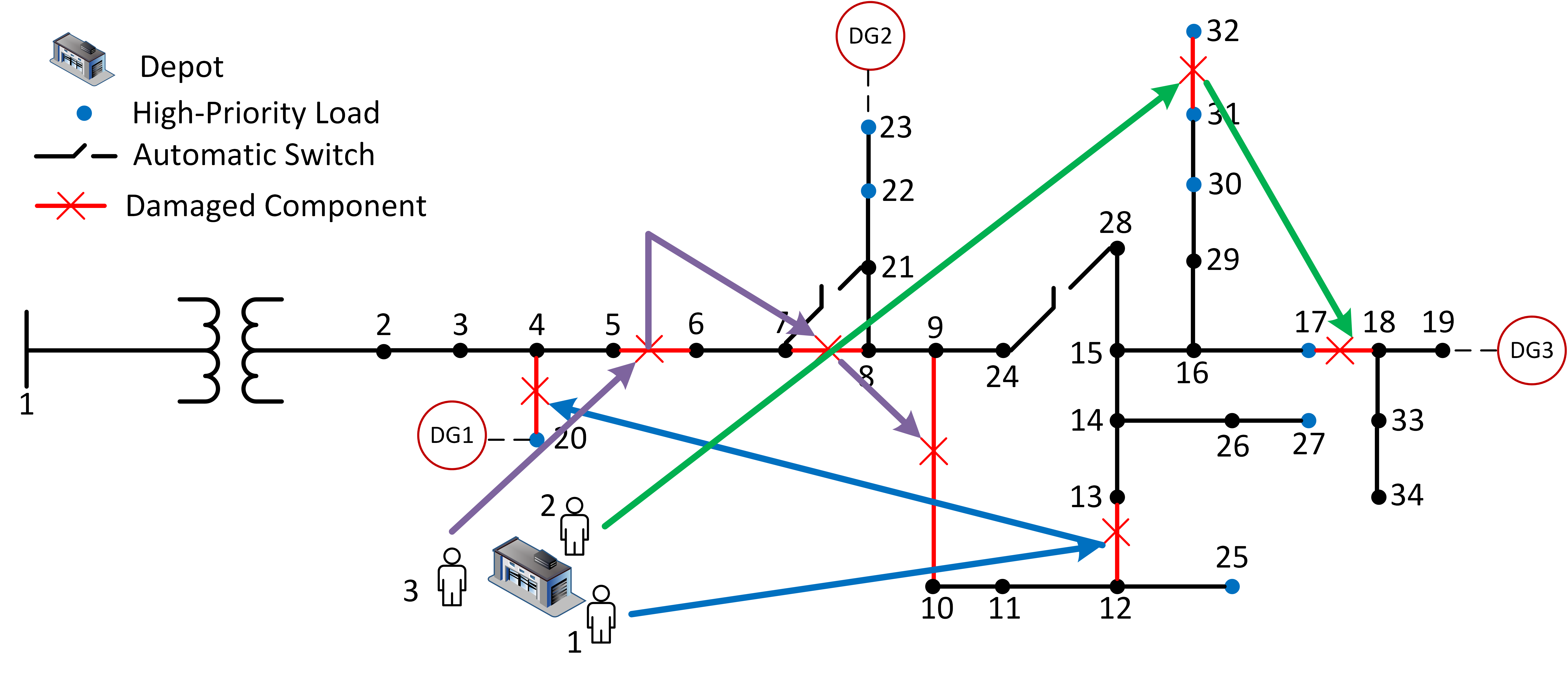}
    \caption{\textcolor{black}{Routing solution for the IEEE 34-bus network obtained by D-PH.}}\label{phsol}
\end{figure}

The aim of this test is to analyze and visualize the D-PH algorithm. Since there are 7 damaged lines and 3 crews, the algorithm requires 3 dispatch cycles. The algorithm converges after 10 minutes, where dispatch cycles 1, 2, and 3 converges after 5, 3, and 2 minutes, respectively. The routing solution is shown in Fig. \ref{phsol}. In the first dispatch cycle, Lines 5-6, 12-13, and 31-32 are selected as critical lines. Repairing line 5-6 provides a path for the power flow coming from the substation. Line 31-32 is prioritized as it is connected to a high-priority load. Line 12-13 is repaired to provide electricity to the lower portion of the network. Line 4-20 is repaired after Line 12-13 as DG1 can provide energy to the load at bus 20 temporarily before the line is repaired.

Next, we present a detailed solution of the second-stage variables for one possible realization, we use Scenario 1 from Table \ref{RTscenarios}. The first-stage solution (crew routing) is shown in Fig. \ref{phsol}, while some of the second-stage variables, including switching operation and DG output, are detailed in Table \ref{DC_T2}. Switch 24-28 is turned on so that DG2 can supply part of the network on the right-hand side. In this scenario, the first line repaired is 31-32, but the load at bus 32 is not served as DG2 is at its limit. Line 12-13 is repaired next and the load at bus 10 is restored. Switch 7-21 remains off until line 5-6 is repaired, to provide a path for the power coming from the substation. The substation restores eight loads at this point (4 AM), while loads at buses 11, 16, and 24 are not restored until the next hour due to the higher demand caused by CLPU. Switch 7-21 and 24-28 are turned off once line 7-8 and line 9-10 are repaired, respectively. Note that by using switches 7-21 and 24-28, all loads are served before repairing lines 7-8 and 9-10. Finally, the back-up DGs are turned off since the loads can be supplied by the substation.

\begin{table}[htbp]
  \centering
  \small
  \caption{\textcolor{black}{Switch status, DG output, and sequence of repairs for the IEEE 34-bus feeder}}
  \vspace{-0.2cm}
    \begin{tabular}{ccccccc}
    \hline
\rule{0pt}{2.2ex}     \multirow{2}[2]{*}{Time} & \multirow{2}[2]{*}{SW 7-21} & \multirow{2}[2]{*}{SW 24-28} & DG1  & DG2  & DG3  & Repaired\\
          &       &       & (kW)  & (kW)  & (kW)  & Component \\
    \hline
\rule{0pt}{2.2ex}     0:00  & 0     & 1     & 74.9  & 143   & 38.7  &  \\
    1:00  & 0     & 1     & 77.5  & 148   & 40    & Line 31-32 \\
    2:00  & 0     & 1     & 67.3  & 149   & 34.8  & Line 12-13 \\
    3:00  & 0     & 1     & 66.4  & 145   & 34.3  & Line 5-6 \\
    4:00  & 1     & 1     & 65.8  & 150     & 34    &  \\
    5:00  & 1     & 1     & 65.8  & 150      & 34    & Line 17-18,4-20 \\
    6:00  & 1     & 1     & 150   & 150   & 150   &  \\
    7:00  & 1     & 1     & 0     & 0     & 0     & Line 7-8 \\
    8:00  & 0     & 1     & 0     & 0     & 0     &  \\
    9:00  & 0     & 1     & 0     & 0     & 0     &  \\
    10:00 & 0     & 1     & 0     & 0     & 0     &  \\
    11:00 & 0     & 1     & 0     & 0     & 0     & Line 9-10 \\
    12:00 & 0     & 0     & 0     & 0     & 0     &  \\
    \hline
    \end{tabular}%
  \label{DC_T2}%
\end{table}%

To show the importance of considering uncertainty in the problem, we calculate the expected value of perfect information (EVPI) and the value of the stochastic solution (VSS). EVPI is the difference between the wait-and-see (WS) and the stochastic solutions. It represents the value of knowing the future with certainty. WS is the expected value of reacting to random variables with perfect foresight. It is obtained by calculating the mean of all deterministic solutions of the scenarios. VSS indicates the benefit of including uncertainty in the optimization problem. VSS is the difference between the stochastic solution and the expected value solution (EEV). To obtain EEV, we first solve the deterministic problem using the expected value (EV) of the random variables, where the average repair time is 4 hours and the load forecast error is zero. Then we set the first-stage variable as a fixed parameter and solve the stochastic problem to find the value of EEV. {\color{black}Furthermore, the expected energy not supplied (EENS) is calculated as follows:

\begin{equation}
\textrm{EENS} = \sum_{\forall s} \textrm{Pr(s)} \big(\sum_{\forall t}\sum_{\forall i}(1-y_{i,t,s})P^D_{i,t,s}\big)
\end{equation}}

\begin{figure}[h!]
\vspace{-0.2cm}
\setlength{\abovecaptionskip}{0pt} 
\setlength{\belowcaptionskip}{0pt} 
  \centering
    \includegraphics[width=0.49\textwidth]{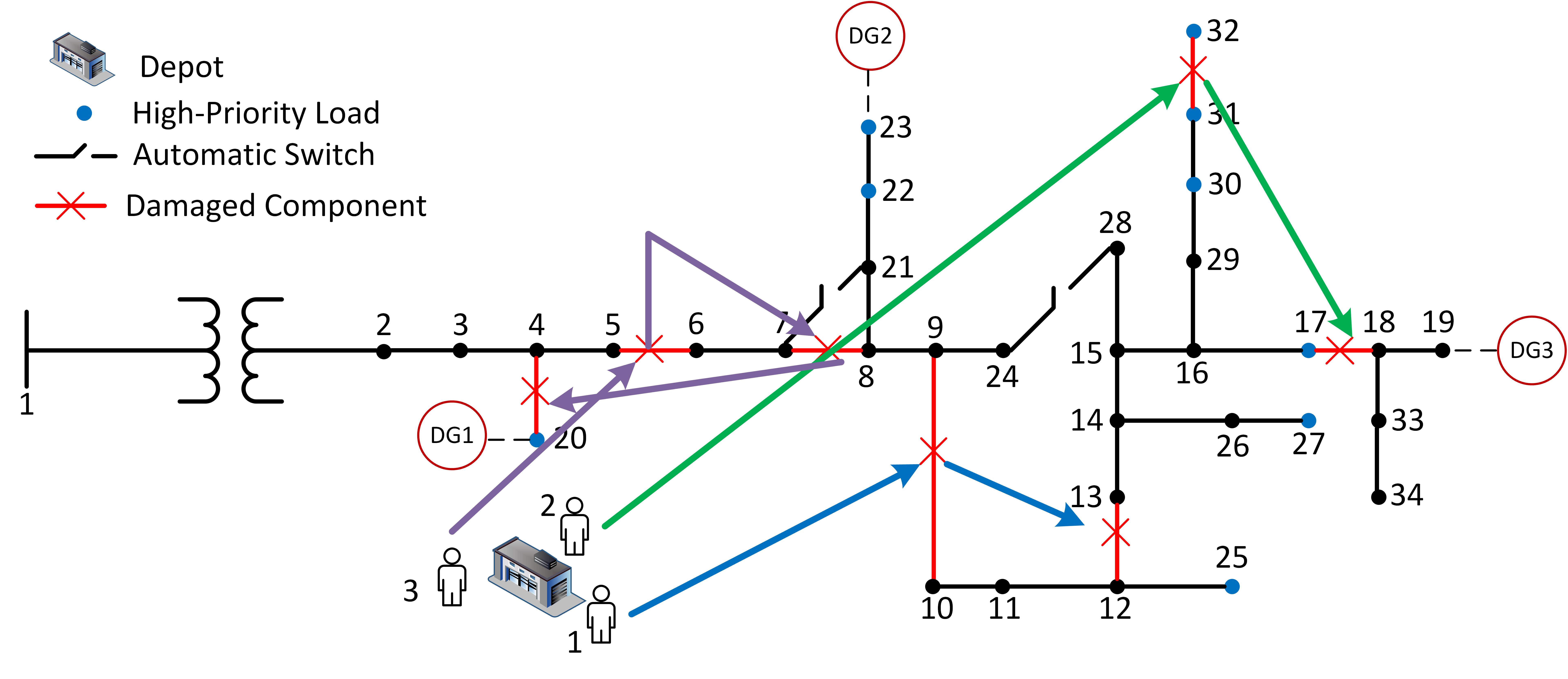}
    \caption{\textcolor{black}{Routing solution obtained by using the expected values.}}\label{EVsol}
\end{figure}

\begin{figure}[h!]
\vspace{-0.2cm}
\setlength{\abovecaptionskip}{0pt} 
\setlength{\belowcaptionskip}{0pt} 
  \centering
    \includegraphics[width=0.49\textwidth]{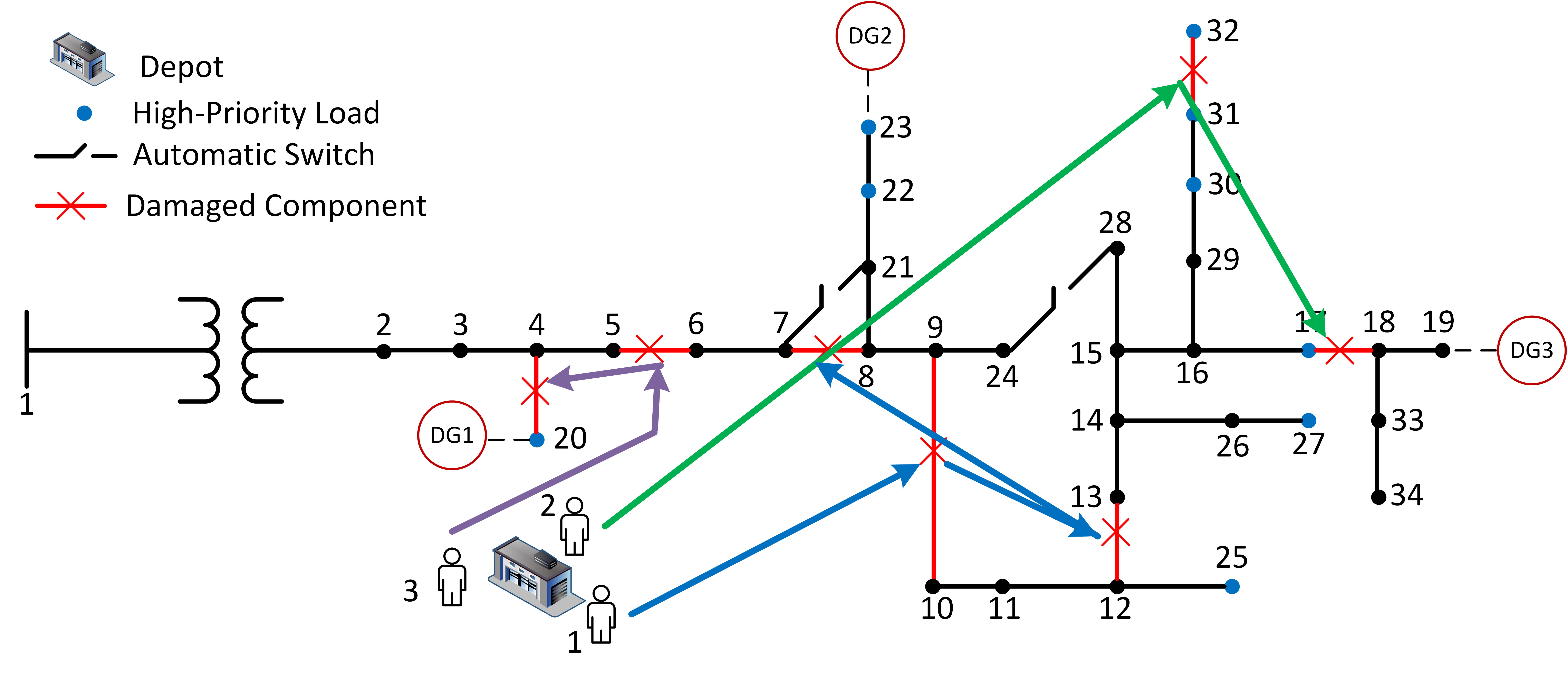}
    \caption{\textcolor{black}{Routing solution obtained by solving the extensive form.}}\label{EFsol}
\end{figure}

The route obtained by solving the deterministic problem with average repair time and zero load forecast error is shown in Fig. \ref{EVsol}. EEV is then found to be 30524.13 and the EENS for this routing plan is 1907.5 kWh, as shown in Table \ref{table34}. By solving the extensive form of the S-DSRRP using Pyomo with CPLEX solver, we obtained the routes shown in Fig. \ref{EFsol}, after 25 hours. Observe that the difference between Fig. \ref{EVsol} and Fig. \ref{EFsol} lies around line 4-20. Repairing line 4-20 early gives DG1 the opportunity to support the substation and meet the higher demand caused by CLPU and the high forecast error. The importance of line 4-20 and DG1 is not captured in the EEV solution as the uncertainty is not considered in the decision making process. D-PH algorithm achieved a solution close to the EF solution in 10 minutes, with EENS 21.2 kWh lower than the one obtained for EF. The relative gap is obtained by comparing the objective of the different methods to the solution obtained using EF, which is only 0.1\% for D-PH. The same route as D-PH is obtained by solving the complete problem (29) using the PH algorithm, but the computation time increases to 27 minutes. 
Though D-PH has a slightly lower objective value than EF, the computation time is improved considerably. {\color{black}Furthermore, the results show the advantage of using PH over EF, as the computation time for EF is 25 hours, whereas PH converges in 27 minutes.}

 \begin{table}[htbp!]
  \centering
  \small
  \caption{Results of the stochastic simulation on the IEEE 34-bus feeder, with 7 damaged components}
  \vspace{-0.2cm}
    \begin{tabular}{ccccccc}
    \hline
\rule{0pt}{2.2ex}          & \textcolor{black}{$\zeta$}& CT & VSS  & EVPI  &  \%Gap & \textcolor{black}{EENS (kWh)}\\
    \hline
\rule{0pt}{2ex}    EEV   & 30524.13 & 257 s & N/A    & N/A    &       0.3\% & 1907.5\\
\rule{0pt}{2ex}    D-PH  & 30588.18 & 10 min & 64.05   &   94.87   &  0.1\% & 1862.0\\
\rule{0pt}{2ex}    PH    & 30588.18 & 27 min & 64.05   &   94.87 &  0.1\% & 1862.0\\
\rule{0pt}{2ex}    EF    & 30617.47 & 25 h & 93.34   &   65.58    &        N/A & 1840.8\\
\rule{0pt}{2ex}    WS    & 30683.05 & 18 min & N/A    & N/A    &       N/A & 1800.4\\
    \hline
    \end{tabular}%
  \begin{tablenotes}
      \small
\item{$\zeta$: objective value (weighted kWh); CT: computation time}
    \end{tablenotes}   
  \label{table34}%
\end{table}%

\subsection{Case II: IEEE 8500-bus distribution feeder}

\begin{figure}[h!]
\vspace{-0.2cm}
\setlength{\abovecaptionskip}{0pt} 
\setlength{\belowcaptionskip}{0pt} 
  \centering
    \includegraphics[width=0.47\textwidth]{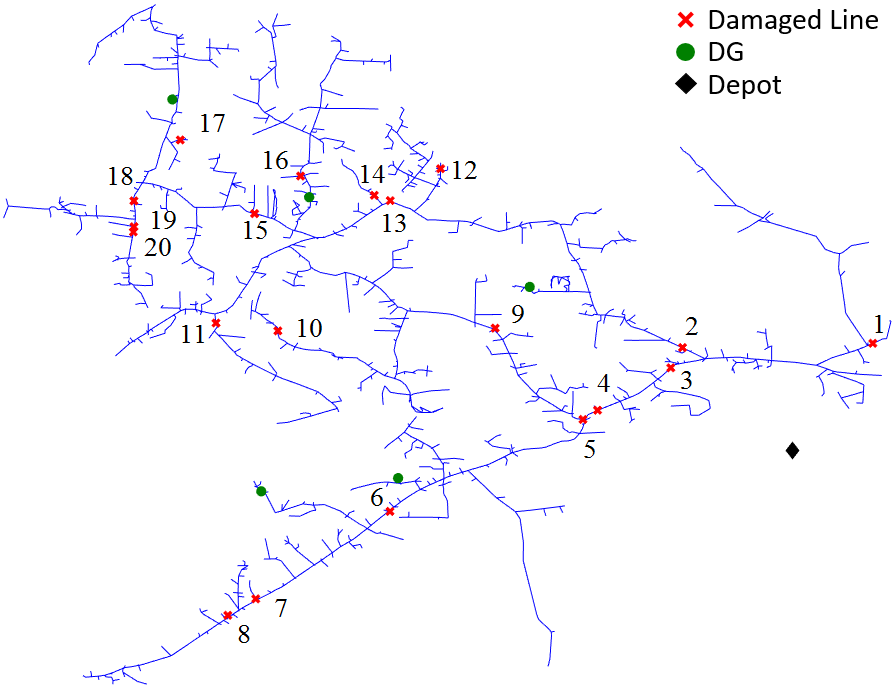}
    \caption{8500-bus IEEE distribution network with 20 damaged lines.}\label{8500fig}
\end{figure}

\begin{figure}[h!]
\vspace{-0.2cm}
\setlength{\abovecaptionskip}{0pt} 
\setlength{\belowcaptionskip}{0pt} 
  \centering
    \includegraphics[width=0.49\textwidth]{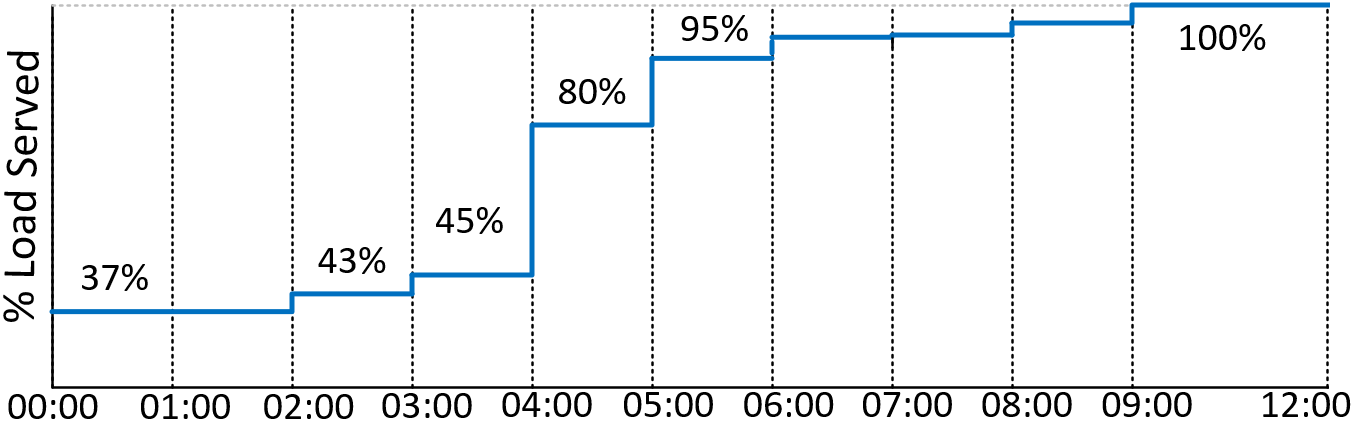}
    \caption{Percentage of load served for the 8500-bus test case.}\label{Load_Rest}
\end{figure}

The IEEE 8500-bus feeder test case, shown in Fig. \ref{8500fig}, is used to examine the scalability of the developed approach for large networks. Five 500 kW DGs are randomly installed in the network. {\color{black}The potential loops in the network are identified using a depth-first search method \cite{Coreman} in MATLAB to form the radiality constraint. There are 5 loops in the network, which are found in 60.72 seconds.} It is assumed that there are 6 crews and 20 arbitrarily selected damaged lines, labeled in Fig. \ref{8500fig}. Monte Carlo sampling is used to generate 1000 random scenarios, which are reduced to 30 using SCENRED2. Since there are 6 crews and 20 damaged lines, the D-PH has four dispatch cycles. The complete routing solution is obtained after 79 minutes, where the 4 dispatch cycles converged after 23, 25, 18, and 13 minutes. {\color{black}The alternative methods, i.e., EEV, EF, and PH, did not converge to a feasible solution after 24 hours}. The routing solution obtained using D-PH is shown in Table \ref{8500table}. Fig. \ref{Load_Rest} shows the change in percentage of load supplied for one sample scenario. By changing the topology of the network and using the backup DGs, 37\% of the loads can be served. The number of served loads start to increase as the crews repair the damaged components, and 95\% of the loads are restored after five hours.


\vspace{-0.2cm}
\begin{table}[h!]
\small
  \centering
  \caption{Routing solution for the 8500-bus test case}
  \vspace{-0.2cm}
    \begin{tabular}{l|ccccccccccc}
    \hline
     \rule{0pt}{2ex}Crew & \multicolumn{10}{c}{Route}\\
    \hline
	\rule{0pt}{2ex}Crew 1 & Depot & $\rightarrow$ & 1 & $\rightarrow$ & 10 & $\rightarrow$ & 9 & $\rightarrow$ & Depot\\
	Crew 2 & Depot & $\rightarrow$ & 15 & $\rightarrow$ & 14 & $\rightarrow$ & 13 & $\rightarrow$ & Depot\\
	Crew 3 & Depot & $\rightarrow$ & 18 & $\rightarrow$ & 7 & $\rightarrow$ & 4 & $\rightarrow$ & Depot\\
	Crew 4 & Depot & $\rightarrow$ & 19 & $\rightarrow$ & 20 & $\rightarrow$ & 6 & $\rightarrow$ & Depot\\
	Crew 5 & Depot & $\rightarrow$ & 11 & $\rightarrow$ & 2 & $\rightarrow$ & 16 & $\rightarrow$ & 12 & $\rightarrow$ & Depot\\
	Crew 6 & Depot & $\rightarrow$ & 5 & $\rightarrow$ & 17 & $\rightarrow$ & 8 & $\rightarrow$ & 3 & $\rightarrow$ & Depot\\
    \hline
    \end{tabular}%
  \label{8500table}%
\end{table}%

\vspace{0.2cm}

{\color{black}
To test whether the scenario set can represent the uncertainties, we apply one of the solution stability tests presented in \cite{Kaut2007}. We perform a sensitivity analysis with different numbers of scenarios for the IEEE 8500-bus system. The stochastic problem is solved to compare the objective values under different numbers of scenarios. The solution is stable if the deviation of these objective values is small \cite{Kaut2007}. The largest number of scenarios we consider is 100. The results are shown in Fig. \ref{stability}. It can be seen that the variation of these objective values is very small, thus, the presented method is stable. This shows that using 30 scenarios can represent the uncertainties in the problem.

\begin{figure}[h!]
\vspace{-0.2cm}
\setlength{\abovecaptionskip}{0pt} 
\setlength{\belowcaptionskip}{0pt} 
  \centering
    \includegraphics[width=0.49\textwidth]{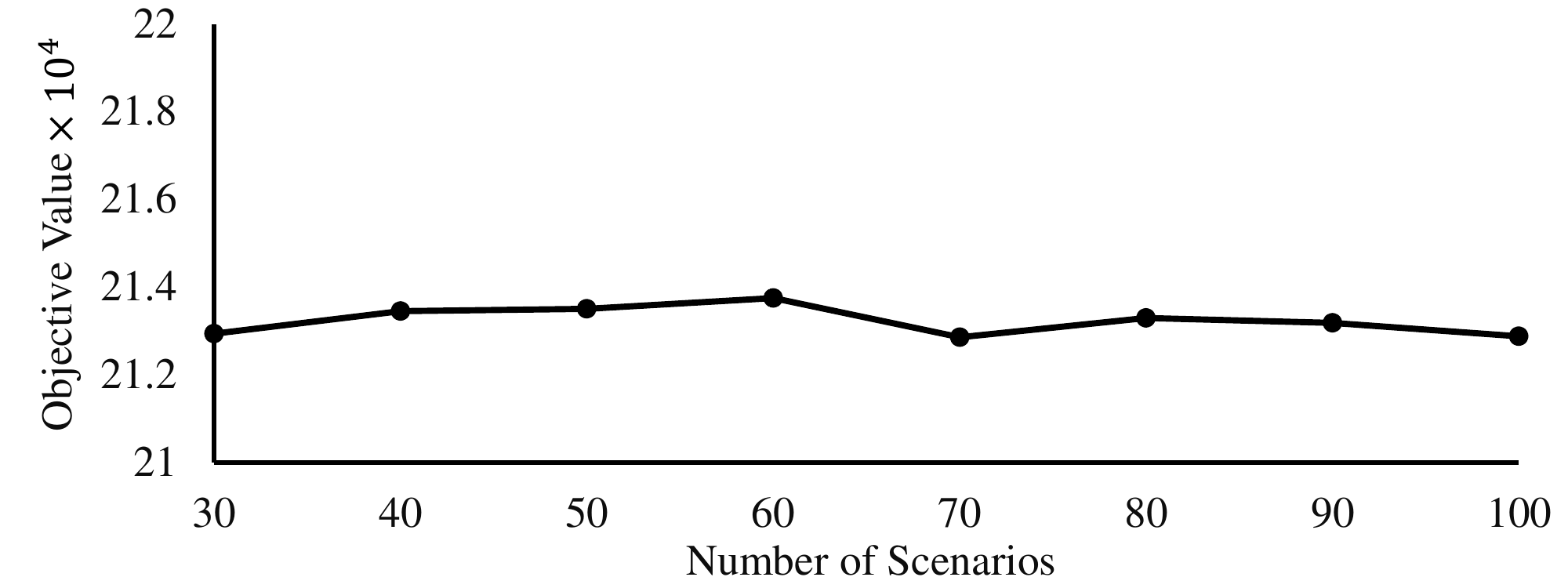}
    \caption{Sensitivity analysis of optimal objective value versus the number of scenarios.}\label{stability}
\end{figure}
}
\section{Conclusion}\label{chap:6}

In this paper, we proposed a two-stage stochastic approach for the repair and restoration of distribution networks. The scenarios are generated using Monte Carlo sampling, considering the uncertainty of the repair time and load. We developed a decomposition approach to solve the stochastic problem. The approach starts with identifying the critical components to repair in its first subproblem, and then routes the crews in the second subproblem. Both subproblems are formulated as two-stage stochastic programs. Parallel Progressive Hedging is employed in the algorithm where the subproblem for each scenario is solved separately. For small cases, the proposed method provides solutions that have similar quality as the one found by solving the extensive form, while the computational burden is significantly reduced. The proposed approach managed to solve large cases in a reasonable time while other methods did not provide a feasible solution within 24 hours. The results demonstrate the effectiveness of the proposed approach in balancing computational burden and solution quality.

\vspace{-0.2cm}

\begin{IEEEbiography}[{\includegraphics[width=1in,height=1.25in,clip,keepaspectratio]{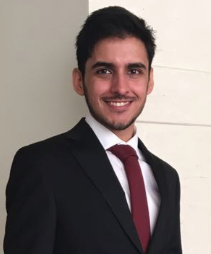}}]{Anmar Arif}
(S'16) is currently pursuing his Ph.D. in the Department of Electrical and Computer Engineering, Iowa State University, Ames, IA. He received his B.S. and Masters degrees in Electrical Engineering from King Saud University and Arizona State University in 2012 and 2015, respectively. Anmar was a Teaching Assistant in King Saud University, and a Research Assistant in Saudi Aramco Chair In Electrical Power, Riyadh, Saudi Arabia, 2013. His current research interest includes power system optimization, outage management, and microgrids.
\end{IEEEbiography}
\vspace{-0.2cm}
\begin{IEEEbiography}[{\includegraphics[width=1in,height=1.25in,clip,keepaspectratio]{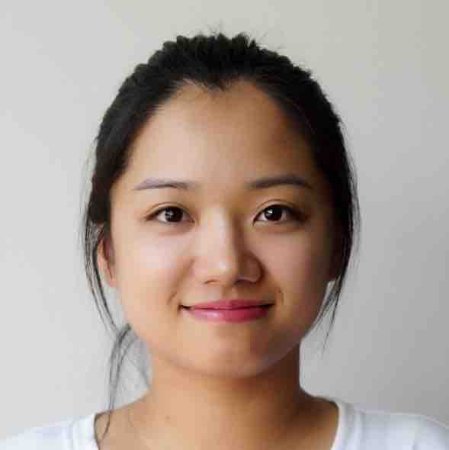}}]{Shanshan Ma}
(S'16) received the B.S. degree in information and electrical engineering from Zhejiang University City College, Hangzhou, China, in 2012, and the M.S. degree from the Department of Electrical Engineering and Computer Science, South Dakota State University, Brookings, SD, USA, in 2015. She is currently pursuing the Ph.D. degree with the Department of Electrical and Computer Engineering, Iowa State University, Ames, IA, USA. Her current research interests include self-healing resilient distribution systems, and microgrids.
\end{IEEEbiography}


\begin{IEEEbiography}[{\includegraphics[width=1in,height=1.25in,clip,keepaspectratio]{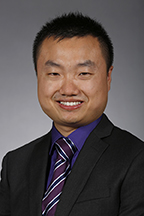}}]{Zhaoyu Wang}
(M'15) received the B.S. and M.S. degrees in electrical engineering from Shanghai Jiaotong University in 2009 and 2012, respectively, and the M.S. and Ph.D. degrees in electrical and computer engineering from the Georgia Institute of Technology in 2012 and 2015, respectively. He is the Harpole-Pentair Assistant Professor with Iowa State University. 
He was a Research Aid with Argonne National Laboratory in 2013, and an Electrical Engineer with Corning Inc. in 2014. His research interests include power distribution systems, microgrids, renewable integration, power system resiliency, demand response, and voltage/VAR control. Dr. Wang was a recipient of the IEEE PES General Meeting Best Paper Award in 2017 and the IEEE Industrial Application Society Prize Paper Award in 2016. He serves as the Secretary of IEEE PES Awards Subcommittee. He is an Editor of the IEEE TRANSACTIONS ON SMART GRID and the IEEE POWER ENGINEERING LETTERS. His research projects are currently funded by the U.S. National Science Foundation, the U.S. Department of Energy, National Laboratories, PSERC, and Iowa Economic Development Agency and Industry.
\end{IEEEbiography}

\begin{IEEEbiography}[{\includegraphics[width=1in,height=1.25in,clip,keepaspectratio]{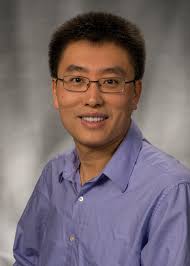}}]{Jianhui Wang}
(M'07-SM'12) received the Ph.D. degree in electrical engineering from the Illinois Institute of Technology, Chicago, IL, USA, in 2007. He is an Associate Professor with Bobby B Lyle School of Engineering, Southern Methodist University, Dallas, TX, USA. He is also the Section Lead of Advanced Power Grid Modeling with the Energy Systems Division, Argonne National Laboratory, Argonne, IL, USA, an Affiliate Professor with Auburn University, Auburn, AL, USA, and an Adjunct Professor with the University of Notre Dame, Notre Dame, IN, USA.
Dr. Wang was a recipient of the IEEE Power and Energy Society (PES) Power System Operation Committee Prize Paper Award in 2015. He is the Secretary of the IEEE PES Power System Operations Committee, an Associate Editor of the Journal of Energy Engineering, an Editorial Board Member of Applied Energy, the Editor-in-Chief of the IEEE Transactions on Smart Grid, and an IEEE PES Distinguished Lecturer.
\end{IEEEbiography}

\begin{IEEEbiography}[{\includegraphics[width=1in,height=1.25in,clip,keepaspectratio]{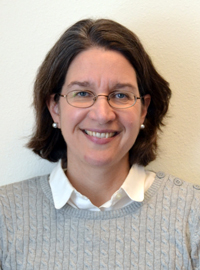}}]{Sarah M. Ryan}
(M'09) received the Ph.D. degree from The University of Michigan, Ann Arbor, MI, USA. She is currently the Joseph Walkup Professor in the Department of Industrial and Manufacturing Systems Engineering at Iowa State University, Ames, IA, USA. Her research applies stochastic modeling and optimization to the planning and operation of service and manufacturing systems.
\end{IEEEbiography}

\begin{IEEEbiography}[{\includegraphics[width=1in,height=1.25in,clip,keepaspectratio]{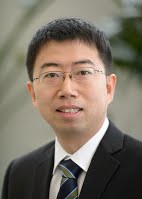}}]{Chen Chen}
(M'13) received the B.S. and M.S. degrees in electrical engineering from Xian Jiaotong University, Xian, China, in 2006 and 2009, respectively, and the Ph.D. degree in electrical engineering from Lehigh University, Bethlehem, PA, USA, in 2013. During 2013-2015, he worked as a Postdoctoral Researcher at the Energy Systems Division, Argonne National Laboratory, Argonne, IL, USA. Dr. Chen is currently a Computational Engineer with the Energy Systems Division at Argonne National Laboratory. His primary research is in optimization, communications and signal processing for smart electric power systems, cyber-physical system modeling for smart grids, and power system resilience.
\end{IEEEbiography}
\enlargethispage{-2.2cm}\vfill

\end{document}